\documentclass[12pt]{amsart}

\usepackage[arrow,matrix]{xy}
\usepackage{amssymb}

\advance\textheight1cm

\advance\headheight1.15pt

\let\:\colon
\let\epsilon\varepsilon
\let\xto\xrightarrow
\let\phi=\varphi
\let\theta=\vartheta
\let\cal=\mathcal
\let\Bbb=\mathbb
\let\frak=\mathfrak
\def\opn#1#2{\def#1{\operatorname{#2}}}
\opn\gp{gp} \opn\Max{Max} \opn\Ker{Ker} \opn\Ext{Ext}
\opn\conv{conv} \opn\chara{char} \opn\n{n} \opn\h{h} \opn\GL{GuL}
\opn\SL{SL} \opn\sn{sn} \opn\int{int} \opn\End{End}
\opn\rank{rank} \opn\Aff{Aff} \opn\Spec{Spec} \opn\Proj{Proj}
\opn\QF{QF} \opn\I{Im} \opn\Hom{Hom} \opn\Aut{Aut} \opn\W{Witt}
\opn\w{w} \opn\int{int} \opn\Pic{Pic}
\def\E{{\cal E}}
\def\BP{{\frak P}}
\def\Pp{{\Bbb P}}
\def\ZZ{{\Bbb Z}}
\def\RR{{\Bbb R}}
\def\NN{{\Bbb N}}
\def\QQ{{\Bbb Q}}
\def\AA{{\Bbb A}}
\def\M{{\frak M}}
\def\N{{\frak N}}
\def\PP{{\cal P}}
\def\O{{\cal O}}

\def\c{\operatorname{\frak c}}

\def\1{^{-1}}
\def\spmatrix#1{\left(\begin{smallmatrix}#1\end{smallmatrix}\right)}
\def\zto{\makebox[0pt]{$0\to\ \ \ \ \ \ $}}
\def\toz{\makebox[0pt]{$\ \ \ \ \ \ \to0$}}

\newtheorem{lemma}{Lemma}[section]
\newtheorem{corollary}[lemma]{Corollary}
\newtheorem{theorem}[lemma]{Theorem}
\newtheorem{proposition}[lemma]{Proposition}

\theoremstyle{definition}
\newtheorem{definition}[lemma]{Definition}
\newtheorem{remark}[lemma]{Remark}

\newtheorem{question}[lemma]{Question}
\newtheorem{conjecture}[lemma]{Conjecture}

\begin{document}

\title{Higher $K$-theory of toric varieties}

\author{Joseph Gubeladze}

\thanks{Supported by the Deutsche Forschungsgemeinschaft, INTAS grant
99-00817 and TMR grant ERB FMRX CT-97-0107}

\address{Department of Mathematics, San Francisco
State University, San Francisco, CA 94132, USA}
\address{A. Razmadze Mathematical Institute, Alexidze St. 1, 380093
Tbilisi, Georgia}

\email{soso@math.sfsu.edu}

\maketitle

\section{Introduction}\label{STATE}

\subsection{Statements}
In \cite{Gu1} we showed that vector bundles on affine toric
varieties are trivial. An easy generalization of this result to
the stable situation yields the equality $K_0(R)=K_0(R[M])$ for
every regular ring $R$ and every {\it seminormal} monoid $M$ (see
\cite{Gu2}, also \cite[Corollary 1.4]{Sw}). The key case is the
class of finitely generated {\it normal} monoids which, up to
isomorphism, can be described as additive submonoids of $\ZZ^r$ of
the type $C\cap\ZZ^r$ for some finite-polyhedral convex rational
cone $C\subset\RR^r$, $r\in\NN$.

The direct higher $K$-theoretic analogue of this fact would be the
equalities $K_i(R)=K_i(R[M])$, $i>0$. But this is harder. In fact,
Srinivas has shown \cite{Sr} that the equality fails even for the
simplest singular toric cone ${\Bbb C}[X^2,XY,Y^2]$ already for
$i=1$. In \cite{Gu4}, using different arguments, we showed that
$K_1(R)\not=K_1(R[M])$ for essentially all finitely generated
submonoids of $\ZZ^r$, not isomorphic to $\ZZ^r_+$.

As a plausible substitute for the equalities $K_i(R)=K_i(R[M])$,
we raised in \cite{Gu2} the following question. Assume $c$ is a
natural number $\geq2$ and $M\subset\ZZ^r$ is a submonoid without
non-trivial units. One has the homothety $c\cdot-:M\to M$,
$m\mapsto cm$. It induces the ring homomorphism $R[-^c]:R[M]\to
R[M]$ and, hence, the group homomorphism $K_i(R[M])\to K_i(R[M])$
which we denote by $c_*$ (for fixed $i$). Now consider an element
$x\in K_i(R[M])$. Is it true that $(c^j)_*(x)\in K_i(R)$ for
$j\gg0$? Speaking loosely: is the multiplicative action of $\NN$
on $K_i(R[M])$ nilpotent?

The conjectural positive answer to this question is strong enough
to contain Quillen's fundamental result on $K$-homotopy invariance
of regular rings (see below) and the aforementioned equality
$K_0(R)=K_0(R[M])$ for {\it arbitrary} seminormal monoid $M$
(Proposition \ref{anders}(a,b)). Moreover, we explain how the
positive answer yields a similar behavior of arbitrary equivariant
closed subsets with respect to the embedded torus (Proposition
\ref{anders}(c)), which should be thought of as a higher version
of Vorst's results \cite{Vo2}. We also show that all of this
generalizes to not necessarily affine toric varieties (Proposition
\ref{global}).

Such a positive answer for $i=1$ is obtained in \cite{Gu2} and for
$i=2$ in \cite{Msh}. In \cite{Gu3} we showed that the answer is
again `yes' for all higher $K$-groups when the cone
$\RR_+M\subset\RR^r$ is {\it simplicial}, i.~e. spanned by
linearly independent vectors. (See also Theorem \ref{simplok}.)
Hence the following conjecture (about a slightly stronger
nilpotence property), which contains all the previous results in a
uniform way:

\begin{conjecture}\label{conj}
Let $M$ be arbitrary commutative, cancellative, torsion free
mo\-no\-id without non-trivial units, $R$ be a regular ring and
$i>0$. Then for every sequence $\c=(c_1,c_2,\dots)$ of natural
numbers $\geq2$ and every element $x\in K_i(R[M])$ one has
$(c_1\cdots c_j)_*(x)\in K_i(R)$ for $j\gg0$ (depending on $x$).
\end{conjecture}

It follows from the geometric approach, developed in \cite{Gu1}
(and adapted here to higher $K$-groups), that this conjecture
admits the following reformulation.

Let $C\subset\RR^r$ be a finite-polyhedral rational cone and
$H\subset\RR^r$ be a rational hyperplane, cutting $C$ into two
non-degenerate subcones $C=C'\cup C''$ in such a way that all
extremal rays of $C$, except one, belong to $C'$ and the
distinguished extremal ray belongs to $C''$. Then Conjecture
\ref{conj} is equivalent to the claim that for an element $x\in
K_i(R[C\cap\ZZ^r])$ and a sequence $\c$ of natural numbers as
above $(c_1\cdots c_j)_*(x)\in\I{\big(}K_i(R[C'\cap\ZZ^r]) \to
K_i(R[C\cap\ZZ^r]){\big)}$ for $j\gg0$ (Lemma \ref{piram}). This
version of Conjecture \ref{conj} could loosely be called
``pyramidal descent''.

The first main result of this paper ``almost'' accomplishes the
pyramidal descent when the coefficient ring is a field $k$
(Theorem \ref{pdescent}). Namely, assuming the induction
hypothesis on $\dim C$, we construct a subring $\Lambda$
(depending on $x$) of the matrix ring
$M_{2\times2}(k[C\cap\ZZ^r])$, whose underlying $k$-vector space
is spanned by a set of monomials (i.~e. matrices with single
non-zero entries from $C\cap\ZZ^r$) containing {\it only one}
element from $C''$ -- all the other monomials in $\Lambda$ come
from $C'$ (notation as above) -- and such that $(c_1\cdots
c_j)_*(x)\in\I{\big(}K_i(R[\Lambda]) \to
K_i(M_{2\times2}(k[C\cap\ZZ^r]))=K_i(R[C\cap\ZZ^r]){\big)}$ for
$j\gg0$.

The ring $\Lambda$ is the ring of endomorphisms of certain rank 2
vector bundle on certain quasiprojective variety. Its construction
involves higher $K$-theoretic analogues of essentially all steps
in \cite{Gu1}, with use of results from \cite{Gu2}. This becomes
possible in combination with Suslin-Wodzicki's solution to the
excision problem and the Mayer-Vietoris sequence for singular
varieties due to Thomason.

If the action of the big Witt vectors on $K_i(A,A^+)$ for graded
not necessarily commutative $k$-algebras $A=A_0\oplus
A_1\oplus\cdots$ satisfies a very natural condition when $\chara
k=0$ (see Question \ref{witt}) then one can also eliminate the
distinguished monomial in $\Lambda$. This would complete the proof
of Conjecture \ref{conj} in 0 characteristic.

Goodwillie's theorem on the relative theories for nilpotent ideals
yields a partial result on Question \ref{witt} in the arithmetic
case which, in combination with Theorem \ref{pdescent} and
Keating's results on congruence-triangular matrices, leads to the
second main result of the paper (Theorem \ref{bp}). Namely, we
present first non-simplicial examples -- actually, an explicit
infinite series of pairwise different combinatorial types of
polyhedral cones, verifying Conjecture \ref{conj} (when $R$ is a
number field) for all higher $K$-groups simultaneously. These are
the cones whose compact cross sections are pyramids or bipyramids
in an iterative way. The simplest nonsimplicial example of such
toric varieties is the cone over the Segre embedding
$\Pp^1\times\Pp^1\to\Pp^3$, the triviality of vector bundles on
which had been a challenge in the 80s \cite{L}\cite[\S13]{Sw}
(before \cite{Gu1}).

\subsection{$K$-homotopy invariance.}
Let $R$ be a ring and let $\NN$ act nilpotently on
$K_i(R[\ZZ_+^r])$, $r\in\NN$. For a natural number $c$ the
endomorphism $R[-^c]:R[\ZZ_+^r]\to R[\ZZ_+^r]$ makes $R[\ZZ_+^r]$
a free module of rank $c^r$ over itself. Thus, using transfer
maps, we see that for any element $x\in K_i(R[\ZZ_+^r])$ there
exists a natural number $j$ such that $c^j\cdot x\in K_i(R)$.
Since the same is true for another natural number $c'$, coprime to
$c$, we get $x\in K_i(R)$.

That there is no way to use the same trick for general finitely
generated monoids is explained by the following fact \cite{BGu}:
for a finitely generated submonoid $M\subset\ZZ^r$ without
non-trivial units, a ring $R$, and a natural number $c>1$ the
endomorphism $R[-^c]:R[M]\to R[M]$ makes $R[M]$ an $R[M]$-module
of finite projective dimension if and only if $M\cong\ZZ_+^s$ for
some $s\in\ZZ_+$.

\subsection{Contents.} The paper is organized as follows.
In \S\ref{GEOM}--\S\ref{POLAR} we make preliminary reductions and
develop convex geometry for monoids (with emphasis on {\it
polarized monoids}). In \S\ref{VARI} we introduce the crucial
variety $X$ and in \S\ref{PDESCENT} we make a reduction in the
study of its $K$-theory. The first main result is proved in
\S\ref{PYDESC}, while \S\ref{WITT}  provides a link with the
actions of Witt vectors. The second main result is proved in
\S\ref{BP}.

\subsection{Preliminaries.}\label{prel}
All the considered monoids $M$ are assumed to be commutative,
cancellative and torsion free. In other words, we assume that the
natural homomorphisms $M\to\gp(M)\to\QQ\otimes\gp(M)$ are
injective, where $\gp(M)$ refers to the universal group associated
to $M$ (the {\it Grothendieck group}, or the {\it group of
differences} of $M$). Equivalently, our monoids can be thought of
as additive submonoids of real spaces ($M\to\RR\otimes\gp(M)$).
This enables us to use polyhedral geometry in their study.

When we treat monoids separately the monoid operation is written
additively. In monoid rings we switch to the multiplicative
notation.

For a monoid $M$ its group of units (i.~e. the maximal subgroup of
$M$) is denoted by $U(M)$. We put $\rank
M=\dim_{\QQ}\QQ\otimes\gp(M)$.

For a subset $W$ of an Euclidean space $\conv(W)$ refers to the
convex hull of $W$.

The rings below, unless specified otherwise, are assumed to be
commutative and with unit. A free rank $n$ module, upon the choice
of a basis, will be thought of as the module of {\it $n$-rows}. In
particular, its elements are multiplied by $n\times n$-matrices
from the right, while the homomorphisms are written on the left.

$\c=(c_1,c_2,\ldots)$ will {\it always} denote arbitrary sequence
of natural numbers $\geq2$.

Finally, $\NN=\{1,2,\ldots\}$, $\ZZ_+=\{0,1,2,\ldots\}$,
$\ZZ_-=\{0,-1,-2,\ldots\}$ and $\RR_+$ refers to the nonnegative
part of $\RR$.

Our general references to toric geometry are \cite{F}\cite{O}. For
a survey of the previous results on $K$-theory of toric varieties
see \cite{Gu5}

\subsection{Acknowledgments.} Thanks are due to Burt Totaro -- the idea of
invoking the Witt vectors' action is his (from his unpublished
notes on Conjecture \ref{conj}), to Mamuka Jibladze for the
numerous discussions, and to Winfried Bruns for inviting me to the
University of Osnabr\"uck and providing excellent working
conditions.

\section{Geometry of monoids, normality and seminormality}\label{GEOM}

Here we recall the relationship between elementary convex geometry
and monoids as developed in \cite{Gu1}\cite{Gu2}. A pure algebraic
alternative is found in \cite[\S4,\S5]{Sw}.

A monoid $M$ is called {\it seminormal} if the following
implication holds
$$
{\big(}m\in\gp(M),\ 2m\in M,\ 3m\in M{\big)}\Rightarrow
{\big(}m\in M{\big)}.
$$
A monoid $M$ is called {\it normal} if $m\in M$ whenever
$m\in\gp(M)$ and $cm\in M$ for some $c\in\NN$.

It is a well known result on monoid domains that for a domain $R$
both the normality and seminormality conditions for $R[M]$ are
equivalent to the corresponding conditions on $R$ and $M$
simultaneously. (A domain $\Lambda$ is called seminormal if the
cancellative and torsion free multiplicative monoid
$\Lambda/U(\Lambda)$ is such.)

Given a finitely generated monoid $M$ with trivial $U(M)$. By
fixing an embedding $M\to\gp(M)\to\RR\otimes\gp(M)$ we can
identify $M$ with an additive submonoid of the Euclidean space
$\RR^{\rank M}$. Clearly, $M$ gets identified with an additive
submonoid of $\QQ^{\rank M}$. We let $C(M)$ denote the cone $M$
spans in $\RR^{\rank M}$, i.~e. $C(M)=\RR_+M$. Then $C(M)$ is a
finite polyhedral, rational, strictly convex (sometimes called
`pointed') cone -- a classical observation in toric geometry. An
element $x\in M$ will be called {\it extremal} if it lies on an
edge (1-dimensional face) of $C(M)$. It is known (see, for
instance, \cite[Theorem 4.5]{Sw}) that there exists a rational
codimension 1 affine subspace $H\subset\RR^{\rank
M}\setminus\{0\}$ such that $C(M)$ is spanned over $O\in\RR^{\rank
M}$ by $\Phi(M)=H\cap C(M)$ -- a finite rational convex polytope.
Of course, $\Phi(M)$ is only defined up to projective
transformation. To any submonoid $M'\subset\QQ_+M$ we associate
the convex subset $\Phi(M')=\RR_+M'\cap H\subset\Phi(M)$. It is
spanned by rational points. In general, when we write $\Phi(L)$
for a monoid $L$ it is assumed that $L\subset\QQ_+M$ for some
finitely generated monoid $L$ with trivial $U(M)$.

A monoid extension $M_1\subset M_2$ is called {\it integral} if a
positive multiple of any element $m\in M_2$ belongs to $M_1$. If
$M_1\subset M_2$ is an integral monoid extension and $M_1$ is
finitely generated, having no nontrivial units, then $M_2$ spans
the same finite polyhedral pointed cone in $\RR^{\rank
M_1}=\RR^{\rank M_2}$. It is also clear that the integrality of
the extension $M_1\subset M_2$ is equivalent to the inclusion
$M_2\subset\QQ_+M_1= \QQ^{\rank M_1}\cap C(M_1)$. The {\it Gordan
Lemma} asserts that

\begin{lemma}\label{Gord}
A finite rank monoid $M$  is finitely generated if and only if
$\gp(M)$ is finitely generated and the cone
$\RR_+M\subset\RR^{\rank M}$ is finite polyhedral rational.
\end{lemma}

Here we do not assume that $U(M)=0$, i.~e. the cone $\RR_+M$ may
not be pointed. A proof in the case when $U(M)=0$ is given in
\cite[Theorem 4.4]{Sw} and the general case reduces to this case
by considering a suitable polyhedral partition of $\RR_+M$.

For a monoid $M$ its {\it normalization} is defined by
$$
\n(M)=\{m\in\gp(M)|\ \exists c_m\in\NN\ \ c_mm\in M\}
$$
and the {\it seminormalization} is defined by
$$
\sn(M)=\{m\in\gp(M)|\ 2m,3m\in M\}.
$$
They satisfy the obvious universality conditions. Observe that
$\n(M)=C(M)\cap\gp(M)$ for a finite rank monoid $M$.

\begin{lemma}\label{filtun}
Both nor\-ma\-li\-za\-ti\-on and se\-mi\-nor\-ma\-li\-za\-ti\-on
pre\-ser\-ve fi\-ni\-te ge\-ne\-ra\-ti\-on. More\-over, any normal
monoid is a filtered union of finitely generated normal monoids
and a seminormal monoid is a filtered union of finitely generated
seminormal monoids.
\end{lemma}

In fact, the group of differences is a filtered union of its
finitely generated subgroups. This reduces the problem to the case
of finitely generated group of differences. Then we approximate
the cone of our monoid by finite rational polyhedral subcones and
Lemma \ref{Gord} applies.

For a finite polyhedral pointed cone $C\subset\RR^r$, $r\in\NN$ we
denote by $\int(C)$ the relative interior of $C$. Similarly, for a
finite convex polytope $\Phi$ its relative interior is denoted by
$\int(\Phi)$. Assume $M$ is a finitely generated monoid with
trivial $U(M)$. Then $\int(M)$ refers to the interior ideal
$\int(C(M))\cap M\subset M$ (``ideal'' here means
$\int(M)+M\subset\int(M)$) and $M_*$ refers to the interior
submonoid $\int(M)\cup\{0\}\subset M$. A pure algebraic definition
is given in \cite[\S5]{Sw}: $M_*=\{m\in M|\ \forall n\in M\ \
\exists c\in\NN\ \ cm-n\in M\}$.

More generally, if $M$ is a monoid such that the cone
$C(M)=\RR_+M\subset\RR^{\rank M}$ is strictly convex finite
polyhedral, we will use the notations $\int(M)=\int(C(M))\cap
M\subset M$ and $M_*=\int(M)\cup\{0\}$.

One more notation. Assume $M$ is a monoid and $W\subset\Phi(M)$ is
a convex subset. Then we put $M(W)=\RR_+W\cap M$. (Here $\RR_+W$
is the cone with vertex at the origin and spanned by $W$.) In
particular, $M_*=M(\int(\Phi(M)))$.

\begin{lemma}\label{int}
Let $M$ be a monoid whose cone $C(M)$ is finite polyhedral and
pointed (hence $U(M)$=0).
\begin{itemize}
\item[(a)]
$\sn(M_*)=\sn(M)_*=\n(M_*)=\n(M)_*$. In particular, if $M$ is
seminormal then $M_*$ is normal.
\item[(b)]
Assume $M$ is finitely generated and normal and $W\subset\Phi(M)$
is a convex subset such that $\dim W=\dim\Phi(M)$ (i.~e. $\dim
W=\rank(M)-1$). Then $M(W)$ contains a free basis of $\gp(M)$.
Equivalently, there is a rational simplex $\Delta\subset W$ such
that $M(\Delta)\approx\ZZ^{\rank M}_+$ and
$\gp(M(\Delta))=\gp(M)$.
\item[(c)]
Assume $W\subset\Phi(M)$ is a convex subset such that $\dim
W=\dim\Phi(M)$. Then $\gp(M(W))=\gp(M)$.
\item[(d)]
A finitely generated monoid $M$ with trivial $U(M)$ can be
embedded into $\ZZ_+^{\rank M}$ in such a way that
$M\subset\ZZ_+^{\rank M}$ induces the equality $\gp(M)=\ZZ^{\rank
M}$.
\end{itemize}
\end{lemma}

\begin{proof}
(a) is proved in \cite[Proposition 1.6]{Gu1}. An alternative proof
(in the finite generation case) is given in \cite[Lemma 6.5]{Sw}.

To see (b), fix arbitrarily a unimodular element $m\in\gp(M)$ such
that $\RR_+m$ meets $W$ in its relative interior. Since $M$ is
normal we have $m\in M$. Complete $m$ to a free basis
$\{m,m_2,\ldots,m_{\rank M}\}$ of $\gp(M)$. If a natural number
$c$ is big enough then all the rays $\RR_+(m_2+cm),\ldots,
\RR_+(m_{\rank M}+cm)$ intersect $W$ in its interior. This is so
because the radial directions of the elements
$m_2+cm,\ldots,m_{\rank M}+cm\in\RR^{\rank M}$ approximate the
direction of $m$ as $c\to\infty$. Therefore, for $c$ big enough we
get the desired basis $\{m,m_2+cm,\ldots,m_{\rank
M}+cm\}\subset\gp(M)$.

The claim (c) follows similarly, using the following observation:
for any element $m\in\int(M(W))$ and $n\in M$ there exists a
natural number $c$ such that $n+cm\in M(W)$, yielding the
inclusion $n=(n+cm)-cm\in\gp(M(W))$.

(d) is proved in \cite[Theorem 1.3.2]{Gu3}. Winfried Bruns
suggested the following alternative short proof based on (b). By
Lemma \ref{filtun} we can additionally assume that $M$ is normal.
Consider the dual cone
$$
C(M)^*=\{\xi\in(\RR\otimes\gp(M))^*|\ \forall x\in C(M)\
\xi(x)\geq0\}
$$
in the dual space
$(\RR\otimes\gp(M))^*=\Hom_{\RR}(\RR\otimes\gp(M)),\RR)$. It is
again a pointed convex polyhedral cone of the same dimension $\dim
C(M)$. By (b) there is a free monoid $F^*\subset C(M)^*$ such that
$\gp(F^*)=\Hom_{\ZZ}(\gp(M),\ZZ)$. Then we have the induced
embedding $M\subset F^{**}$ where
$F^{**}=C(F^*)^*\cap\gp(M)\approx\ZZ_+^{\rank M}$.
\end{proof}

\begin{lemma}\label{inv}
Assume $M$ is a finitely generated normal monoid with trivial
$U(M)$ and let $E\subset C(M)$ be an edge. Assume $t$ is the
(unique) generator of $E\cap M\approx\ZZ_+$. Then
$\ZZ_+(-t)+M\approx\ZZ\times N$ for some finitely generated normal
monoid $N$ for which $\rank N=\rank M-1$ and $U(N)=0$.
\end{lemma}

(See, for instance, \cite[Theorem 1.8]{Gu1}\cite[Lemma 8.7]{Sw}.)

\begin{corollary}\label{pyrapp}
Let $M$ and $t$ be as in Lemma \ref{inv}. Then the submonoid
$(\ZZ_+(-t)+M)\setminus\NN(-t)\subset\ZZ_+(-t)+M$ is a filtered
union of monoids of the type $\ZZ t\times N_c$, $c\in\NN$ where
the $N_c$ are finitely generated normal monoids and $U(N_c)=0$.
\end{corollary}

\begin{proof}
Fix any representation $\ZZ_+(-t)+M\approx\ZZ\times N$. By Lemma
\ref{int}(d) there is a basis $\{b_1,\ldots,b_{\rank N}\}\subset
N$ of the group $\gp(N)$ such that
$\Phi(N)\subset\conv\{b_1,\ldots,b_{\rank N}\}$. For a natural
number $c$ consider the monoid $$
N_c={\big(}\ZZ(b_1-ct)+\cdots+\ZZ(b_{\rank N}-ct){\big)}\cap
{\big(}\ZZ_+(-t)+M{\big)}, $$ the intersection being considered in
$\gp(M)$. Using the fact that the radial directions of the
elements $b_1-ct,\ldots,b_{\rank N}-ct\in\RR^{\rank M}$
approximate that of $-t$ when $c\to\infty$ we obtain the desired
filtered union representation $$ \bigcup_{c\in\NN}(\ZZ
t+N_c)=(\ZZ_+(-t)+M)\setminus\NN(-t). $$
\end{proof}

Below, when the monoid operation will be written multiplicatively,
we will use the notation $t\1M$ for $\ZZ_+(-t)+M$.

\section{Excision}\label{EXCI}

Recall that a non-unital not necessarily commutative ring $I$
satisfies {\it excision} in algebraic $K$-theory if the natural
homomorphisms $K_i(\tilde I,I)\to K_i(A,I)$, $i\in\ZZ$ are
isomorphisms for any ring $A$, containing $I$ as a two-sided
ideal. Here $\tilde I=\ZZ\ltimes I$ refers to the universal ring
obtained by adjoining unit (for details see \cite{SuW}). Actually
we only need to require the excision for positive indices $i$
because for $i\leq0$ it is always satisfied \cite[Theorem
XII.8.3]{Ba}. For such a non-unital ring $I$ and a ring $A$,
containing $I$ as a two-sided ideal, we have the long exact
sequence $$ \cdots\to K_{i+1}(A/I)\to K_i(I)\to K_i(A)\to
K_i(A/I)\to\cdots $$ ($i\in\ZZ$), where $K_*(I)=K_*(\tilde I,I)$.

\begin{theorem}[Suslin-Wodzicki \cite{SuW}]\label{SuW}
A non-unital (not necessarily commutative) ring $I$ satisfies
excision if for any finite system $a_1,\ldots,a_n\in I$ there
exist elements $b_1$,$\ldots$,$b_n$,$u$,$v\in I$ such that
$a_i=b_iuv$ for $i\in[1,n]$ and the left annihilators of $u$ and
$uv$ coincide.
\end{theorem}

In this section we show that ``interior ideals'' of certain monoid
rings satisfy excision in $K$-theory.

For a monoid $M$ we denote by $M^{\c}$ the inductive limit of the
diagram
$$
M\xto{c_1\cdot}M\xto{c_2\cdot}\cdots.
$$
We will think of $M^{\c}$ as the filtered union of the ascending
chain
$$
\frac{M}{c_1}\subset\frac{M}{c_1c_2}\subset\cdots,
$$
the union being considered in $\QQ\otimes\gp(M)$.

\begin{lemma}\label{snrm}
$M^{\c}$ is seminormal.
\end{lemma}

\begin{proof}
Assume $2m,3m\in M^{\c}$ for some $m\in\gp(M^{\c})$. Since 2 and 3
generate the additive monoid $\ZZ_+\setminus\{1\}$ there exists
$j\in\NN$ such that
$$
c_{j+1}m\in\frac M{c_1\cdots c_j}.
$$
But then
$$
m\in\frac M{c_1\cdots c_jc_{j+1}}\subset M^{\c}.
$$
\end{proof}

\begin{lemma}\label{exc}
Let $R$ be a ring and $M$ be a finitely generated monoid with
trivial $U(M)$. Then the ideal $J=R\int(M^{\c})\subset R[M^{\c}]$
satisfies excision in $K$-theory.
\end{lemma}

\begin{proof}
Assume $a_1,\ldots,a_n\in J$ for some $n\in\NN$. Fix any element
$m\in\int(M)$. There exists $j\in\NN$ such that the element $$
m_j=m^{(c_1\cdots c_j)\1}\in M^{(c_1\cdots c_j)\1} $$ satisfies
the condition (in the multiplicative notation):
$$
a_1(m_j)\1,\ldots,a_n(m_j)\1\in R[\int(C(M)) \cap\gp(M^{\c})].
$$
By Lemmas \ref{int}(a) and \ref{snrm} we have
$\int(C(M))\cap\gp(M^{\c})= \int(M^{\c})$. Now Theorem \ref{SuW}
applies because of the equalities in $\int(M^{\c})$: $$
m_j=(m_j)^{c_{j+1}\1}\cdot
{\big(}(m_j)^{c_{j+1}\1}{\big)}^{c_{j+1}-1}. $$
\end{proof}

Based on Lemma \ref{exc} we can make the first reduction in
Conjecture \ref{conj}.

\begin{lemma}\label{inter}
To prove Conjecture \ref{conj} for a regular coefficient ring $R$
it suffices to show that $K_i(R)=K_i(R[M^{\c}_*])$ $(i\in\NN)$ for
all finitely generated normal monoids $M$ with trivial groups of
units.
\end{lemma}

(Clearly, $(M^{\c})_*=(M_*)^{\c}$ and, therefore, we can use the
notation $M^{\c}_*$.)

\begin{proof}
Since any monoid is a filtered union of finitely generated monoids
and $K$-functors commute with filtered inductive limits Conjecture
\ref{conj} reduces to the case of finitely generated monoids. We
want to show that $K_i(R)=K_i(R[L^{\c}])$ assuming the equalities
as in the lemma, where $L$ is any finitely generated monoid with
trivial $U(L)$. By Lemmas \ref{int}(a) and \ref{exc} (the latter
being applied to the extension $R\int(L^{\c})\subset R[L^{\c}_*]$)
one has $K_i(R[L^{\c}_*],R\int(L^{\c}))=0$. By Lemma \ref{exc} the
extension $R\int(L^{\c})\subset R[L^{\c}]$ yields the natural
isomorphisms $K_i(R[L^{\c}])=K_i(R[L^{\c}]/R\int(L^{\c}))$. Put
$A_{-1}=R[L^{\c}]$ and $A_0=R[L^{\c}]/R\int(L^{\c})$.

Let $F\subset\Phi(L)$ be any {\it facet} (codimension 1 face).
Then the ideal $R\int(L^{\c}(F))\subset R[L^{\c}(F)]$ is an ideal
of the bigger ring $A_0$ as well. Using the same arguments as
above we arrive at the equalities $K_i(A_0)=K_i(A_1)$ where
$A_1=A_0/R\int(L^{\c}(F))$. Next we do the same reduction with
respect to another facet of the polytope $\Phi(L)$, and so on. By
the same token we get a sequence of rings $A_0,A_1\ldots,A_s$ such
that
$$
A_k=A_{k-1}/R\int(L^{\c}(F_k))\quad\text{and}\quad
K_i(A_{k-1})=K_i(A_k)
$$
for the corresponding enumeration of the facets
$F_k\subset\Phi(L)$, $k\in[1,p]$ where $F_1=F$.

Thereafter we handle the codimension 2 faces of $\Phi(L)$, and so
on. Finally, by annihilating all non-trivial monomials, we shall
descend to the coefficient ring $R$. (When we treat
$0$-dimensional faces, i.~e. the vertices of $\Phi(L)$, we follow
the convention -- an interior of a point is the point itself.)
Finally, we obtain a sequence of rings $A_{-1},A_0\ldots,A_q=R$
($q\geq p$) such that
$K_i(A_{-1})=K_i(A_0)=\cdots=K_i(A_q)=K_i(R)$.
\end{proof}

\begin{proposition}\label{anders}
Let $R$ be a ring and $M$ a be finitely generated monoid with
trivial $U(M)$.
\begin{itemize}
\item[(a)]
If $\NN$ acts nilpotently on $SK_0(R[M])$ then
$SK_0(R)=SK_0(R[M])$
\item[(b)]
If $\NN$ acts nilpotently on $K_0(R[M])$, $R$ is a domain and both
$R$ and $M$ are seminormal then $K_0(R)=K_0(R[M])$.
\item[(c)]
If $I\subset M$ is a proper ideal (i.~e. $I\not=M$ and $IM\subset
I$) and $\NN$ acts nilpotently on $K_i(R[M])$ ($i\in\NN$) then
$\NN$ acts nilpotently also on $K_i(R[M]/I)$.
\end{itemize}
\end{proposition}

(The action of $\NN$ on $SK_0(R[M])$ and $K_i(R[M]/I)$ and the
corresponding nilpotence is understood in the obvious sense.)

Actually one can drop the finite generation condition on $M$ but
we skip this. The arguments below follow closely
\cite[\S3.3]{Gu3}. We give a shortened version, in which Lemma
\ref{int} is used implicitly several times.

\begin{proof}
(a) It follows easily from the results in \cite{I} that for a ring
extension $A\subset A[b]=B$, such that $b^2,b^3 \in A$ the
homomorphism $SK_0(A)\to SK_0(B)$ is injective.

Let $c$ be a natural number. Consider the endomorphism $-^c:M\to
M$. We fix a basis $\{m_1,\ldots,m_r\}$ of $\gp(M)$ which is a
subset of $\int(M)$. Then $R[\I(-^c)+\ZZ_+m_1]$ is a free
$R[\I(-^c)]$-module of rank $c$. In particular, one easily
concludes (using transfer maps for $K_0$)  that any element $x\in
SK_0(R[M])$, which is mapped to zero in
$SK_0(R[\I(-^c)+\ZZ_+m_1])$, is of $c$-torsion. By the infinite
iteration, involving periodically all the elements
$m_1,\ldots,m_r$, we conclude that every element of
$$
\Ker{\big(}SK_0(R[M])\to
SK_0(R[M+(\ZZ_+m_1+\cdots+\ZZ_+m_r)^{\c'}]){\big)},
$$
where ${\c'}=(c,c,\ldots)$, is annihilated by some positive power
of $c$.

In view of the equality
$\sn([M+(\ZZ_+m_1+\cdots+\ZZ_+m_r)^{\c'}]_*) =M^{\c'}_*$ and the
aforementioned corollary of Ischebeck's results, the `interior
excision' arguments (as in the proof of Lemma \ref{inter}) for
$SK_0$-groups show the implication
\begin{itemize}
\item[]
($\NN$ acts nilpotently on $SK_0(R[M])$) $\Rightarrow$
($SK_0(R[M])/SK_0(R)$ is of $c$-torsion).
\end{itemize}
We are done because the same is true for another natural number
$c'$, coprime to $c$.

(b) Since for seminormal monoids $M$ it is well known that
$\Pic(R)=\Pic(R[M])$ provided $R$ is a seminormal domain, the
claim follows from (a).

(c) One only needs to observe that the limit ring
$$\lim_{j\to\infty}{\big(}R[M]/RI\xto{R[-^{c_j}]/RI}R[M]/RI{\big)}$$
is a ``face'' ring of the type we had in the proof of Lemma
\ref{exc}. This is so because the limit contains no non-zero
nilpotent elements, corresponding to monomials from $M$.
Therefore, the similar arguments as in that proof apply.
\end{proof}

\section{Localizations}\label{LOCAL}

A {\it Karoubi square} is a commutative square of (not necessarily
commutative) rings of the type
$$
\xymatrix{
A\ar[r]\ar[d]&B\ar[d]\\
S\1A\ar[r]&S\1B }
$$
where $S$ is a central multiplicative subset of $A$ which is
regular on $A$ and $B$ and such that $A/sA\to B/sB$ is an
isomorphism for every $s\in S$\cite[\S2]{Sw}. In other words we
require that $A\to B$ is an {\it analytic isomorphism along} $S$
\cite{Vo1}\cite{W1}. The basic fact on Karoubi squares is that it
implies the long exact sequence
$$
\cdots\to K_i(A)\to K_i(S\1A)\oplus K_i(B)\to K_i(S\1B)\to
K_{i-1}(A)\to\cdots
$$
(\cite{Vo1}.) This is shown by comparing the two localization
sequences corresponding to $A\to S\1A$ and $B\to S\1B$, with use
of the equivalence ${\bf H}_S(A)\approx{\bf H}_S(B)$ due to
Karoubi \cite{Ka} (see also \cite[\S1]{W1}).

We will need the following lemma, based on an observation of Swan
\cite[\S10]{Sw}.

Assume $R=(R,\mu,k)$ is a local ring and $M$ is a monoid, not
necessarily finitely generated, for which $\Phi(M)$ is defined.
Let $W\subset\Phi(M)$ be a convex {\it open} subset for which
$\dim W=\dim\Phi(M)$. Put $\M=\mu+R(M(W)\setminus\{1\})$ ($\in\max
R[M(W)]$).

\begin{lemma}\label{kar}
Let $A\subset B$ be two (not necessarily commutative)
$R[M(W)]$-algebras ($R[M(W)]$ is in their centers) such that
$M(W)^{-1}A=M(W)^{-1}B$. Then the commutative square of rings
$$
\xymatrix{
A\ar[r]\ar[d]&B\ar[d]\\
A_\M\ar[r]&B_\M }
$$
is a Karoubi square.
\end{lemma}

\begin{proof}
Consider a maximal ideal ${\frak N}\in\max R[M(W)]$. We claim that
there are only two possibilities: either ${\frak N}=\M$ or ${\frak
N}\cap M(W)=\emptyset$. In fact, if $m\in{\frak N}\cap M(W)$ then
for any $n\in M(W)\setminus\{1\}$ there exists a natural number
$c$ such that $n^cm\1\in M(W)$ (as in the proof of Lemma
\ref{int}(c)). Therefore, $n^c\in\N$ which yields $n\in\N$. Since
the only maximal ideal of $R[M(W)]$ containing
$M(W)\setminus\{1\}$ is $\M$ the claim follows.

Pick $s\in S= R[M(W)]\setminus\M$. We want to show that $A/(s)\to
B/(s)$ is an isomorphism. This can be checked locally on $\max
R[M(W)]$. At $\M$ both the source and the target localize to $0$.
If ${\frak N}\in\max R[M(W)]$ is another maximal ideal then we
know that $M(W)\subset S$. Therefore,
$(A/(s))_{\M}\to(B/(s))_{\M}$ is a further localization of the
identity mapping $M(W)^{-1}A/(s)\to M(W)^{-1}B/(s)$.
\end{proof}

One more localization result.

\begin{lemma}\label{patch}
To prove Con\-jec\-tu\-re \ref{conj} it suf\-fi\-ci\-es to show
that $K_i(R)$ $=$ $K_i(R[M^{\c}_*])$ $(i\in\NN)$ for every local
regular ring $R$ and every finitely generated normal monoid $M$
with trivial $U(M)$.
\end{lemma}

The proof below is based on the following local-global patching
for $K$-groups due to Vorst \cite{Vo1}. It is a higher (stable)
analogue of Quillen's original local-global patching for
projective modules \cite{Q2}.

\begin{theorem}\label{lgp}
Let $A$ be a ring and $x\in K_i(A[T])$ ($T$ a variable). Then
$x\in K_i(A)$ iff $x_\mu\in K_i(A_\mu)$ for any $\mu\in\max A$.
\end{theorem}

{\it Proof of Lemma \ref{patch}.} By Lemma \ref{inter} all we need
is to show that the coefficient ring in the equalities
$K_i(R)=K_i(R[M^{\c}_*])$ admits globalization.

{\it Step 1}. Assume $R$ is any regular ring and $M$ is as in the
lemma. Assume further we have shown that
$K_i(\Lambda)=K_i(\Lambda[L^{\c}_*])$ for any regular local ring
$\Lambda$ and any finitely generated normal monoif $L$ with
trivial $U(L)$. Then by Lemma \ref{inter} we have
$K_i(R_\mu)=K_i(R_\mu[M^{\c}_*\times\ZZ_+^{\c}])$ for all
$\mu\in\max R$. Fix an element $x\in K_i(R[M^{\c}_*\times\ZZ_+])$.
There are elements $a_1,\ldots,a_p\in R$ and a natural number
$j\in\NN$ such that $(c_1\cdots c_j)_*(x)_{a_u}\in K_i(R_{a_u})$,
$u\in[1,p]$ and $a_1R+\cdots+a_pR=R$, where $c_*$ refers to the
endomorphism of the group $K_i(R[M^{\c}_*][\ZZ_+])$ induced by
$-^c:\ZZ_+\to\ZZ_+$. It follows that $(c_1\cdots c_j)_*(x)_{\M}\in
K_i(R[M^{\c}_*]_{\M})$ for any $\M\in\max R[M^{\c}_*]$. By Theorem
\ref{lgp} we conclude $(c_1\cdots c_j)_*(x)\in K_i(R[M^{\c}_*])$.
Next we derive the equality
\begin{equation}
\tag{$*$}K_i(R[M^{\c}_*])=K_i(R[M^{\c}_*][\ZZ_+^{\c}]).
\end{equation}
Pick an element $y\in K_i(R[M^{\c}_*\times\ZZ_+^{\c}])$. It admits
a representation $y=(c_1\cdots c_{j_1})_*(z)$ for some $j_1\in\NN$
and $z\in K_i(R[M^{\c}_*][\ZZ_+])$. By the remarks above, applied
to the sequence $\c'=(c_{j_1+1},c_{j_1+2},\ldots)$, there exists
an index $j_2>j_1$ such that
$$
(c_{j_1+1}\cdots c_{j_2})_*(z)\in K_i(R[M^{\c'}_*]),
$$
which implies $(*)$.

{\it Step 2.}  By Lemma \ref{int}(d) we know that $M$ can be
embedded into a free monoid. As a result, there is a
$\ZZ_+$-grading $R[M]=R\oplus R_1\oplus\cdots$ for which $M$
consists of homogenous elements. It follows that there is a
$\ZZ_+^{\c}$ grading
$$
R[M^{\c}_*]=\bigoplus_{\ZZ_+^{\c}}R_\alpha,\ \ R_0=R.
$$
By Lemma 3.5 in \cite{Gu2} (a straightforward generalization from
polynomial algebras to monoid algebras of s.~c. {\it Weibel
homotopy trick}, see below) for any functor ${\frak F}$ from the
category of rings to that of abelian groups the following
implications holds
$$
{\big(}{\frak F}(R[M^{\c}_*])={\frak
F}(R[M^{\c}_*][\ZZ_+^{\c}]){\big)} \Rightarrow{\big(}{\frak
F}(R)={\frak F}(R[M^{\c}_*]){\big)}.
$$
By $(*)$ we are done. \qed

Recall, that the just mentioned `Weibel homotopy trick' is the
following assertion: for any functor ${\frak F}$ from rings to
abelian groups and any graded ring $A=A_0\oplus A_1\oplus\cdots$
the following implication holds
$$
{\frak F}(A)={\frak F}(A[\ZZ_+])\ \Rightarrow\ {\frak
F}(A_0)={\frak F}(A).
$$

In the proof of Theorem \ref{pdescent} below we will make an
essential use of the Mayer-Vietoris sequence for singular
varieties due to Thomason \cite[Theorem 8.1]{TT}:

\begin{theorem}\label{MV}
Let $U$, $V$ be two quasi-compact open subschemes of a
qua\-si-se\-pa\-ra\-ted scheme. Then there is a natural long exact
Mayer-Vietoris sequence
$$
\cdots\to K_i(U\cup V)\to K_i(U)\oplus K_i(V)\to K_i(U\cap V)\to
K_{i-1}(U\cup V)\to\cdots
$$
\end{theorem}

\begin{remark}\label{quilwald}
Actually, the $K$-groups mentioned here are those of Waldhausen,
associated to the appropriate categories of {\it perfect
complexes} on the given schemes. However, these $K$-groups
coincide with Quillen's $K$-groups for quasiprojective (in
particular, affine) schemes over a commutative ring
\cite[\S3]{TT}. In particular, this is so for schemes of the type
$X$ from \S\ref{VARI} (see Remark \ref{toric}(a)). Accordingly,
all the $K$-groups that appear in this paper are {\it a priori}
assumed to be those of Quillen \cite{Q1}.
\end{remark}

\begin{lemma}\label{monic}
For a ring $R$ and a monic polynomial $f\in R[T]$ the natural
homomorphism $K_i(R[T])\to K_i(R[T,f\1])$ is injective.
\end{lemma}

For $i=0,1$ this is proved in \cite[\S1]{MuP}. We use similar
ideas.

\begin{proof}
First we observe that for any ring $B$ and its two comaximal
elements $b_1,b_2$ Theorem \ref{MV} implies the long exact
sequence
$$
\cdots\to K_i(B)\to K_i(B_{b_1})\oplus K_i(B_{b_2})\to
K_i(B_{b_1b_2})\to K_{i-1}(B)\to\cdots
$$

To prove the lemma it is enough to show that $K_i(R[T])\to
K_i(R[T,T\1f\1])$ is injective.

Let $f=T^n+a_{n-1}T^{n-1}+\cdots+a_0$. We write $f=gY^{-n}$, where
$Y=T\1$ and $g=1+a_{n-1}Y+\cdots+a_0Y^n$. We have the equality
$R[T,T\1f\1]=R[Y,Y\1g\1]$. Since $Y$ and $g$ are comaximal in
$R[Y]$ the observation above yields the exact sequence $$
K_i(R[Y])\to K_i(R[Y^{\pm1}])\oplus K_i(R[Y,g\1])\to
K_i(R[Y,Y\1g\1]). $$ By the {\it Fundamental Theorem} \cite{Gra}
we have the embeddings $$ K_i(R[T])\to K_i(R[T^{\pm1}]),\ \
K_i(R[Y])\to K_i(R[Y^{\pm1}])= K_i(R[T^{\pm1}]) $$ with $K_i(R)$
the intersection of their images. Thus, the exact sequence above
implies $$ \Ker(K_i(R[T])\to K_i(R[T,T\1f\1]))=\Ker(K_i(R)\to
K_i(R[T,T\1f\1])). $$ Now the proof is completed by Lemma 1.1 in
\cite{MuP} which says that for any functor $\frak F$ {\it with
transfers} from rings to abelian groups the homomorphisms of the
type ${\frak F}(R)\to {\frak F}(R[T,h\1])$, $h$ monic, are always
injective.
\end{proof}

Next we show that Conjecture \ref{conj} admits a natural
generalization to not necessarily affine (normal) toric varieties.

Let $X$ be a (normal) toric variety over a regular ring $R$. The
multiplicative actions of $\NN$ on the affine toric varieties in
the standard open covering are compatible. In particular, we have
the actions $\NN\times X\to X$ and $\NN\times(X\times\AA^1_R)\to
X\times\AA^1_R$, where we mean the trivial action on $\AA^1_R$.
Denote by $K_i(X)^{\c}$ and $K_i(X\times\AA^1_R)^{\c}$ the
inductive limits of the corresponding $K_i$-groups with respect to
the successive endomorphisms, given by $\c$.

\begin{proposition}\label{global}
Conjecture \ref{conj} is equivalent to the equalities
$K_i(X\times\AA^1_R)^{\c}=K_i(X)^{\c}$ for all quasiprojective
toric varieties $X$.
\end{proposition}

\begin{proof}
If $X$ is affine then $X=Y\times T$ for some affine toric variety
$Y$, whose polyhedral cone is strictly convex, and some torus $T$.
Using the Fundamental Theorem, the desired equality is easily
derived in the affine case from Conjecture \ref{conj}.

Next we reduce the general case to the affine case as follows. Let
$X$ be given by a fan $\mathcal F$. Pick a maximal cone
$\sigma\in\cal F$ and let ${\cal F}-\sigma$ denote the fan
determined by the remaining maximal cones. Denote by $U$ and $V$
the open toric subvarieties of $X$, corresponding to $\sigma$ and
${\cal F}-\sigma$. Theorem \ref{MV} yields the long exact sequence
$$\cdots\to K_{i+1}(U\cap V)\to K_i(X)\to K_i(U)\oplus K_i(V)\to
K_i(U\cap V)\to\cdots $$

By passing to the inductive limits we get the exact sequence (in
the self explanatory notation)
$$\cdots\to K_{i+1}(U\cap V)^{\c}\to K_i(X)^{\c}\to K_i(U)^{\c}\oplus
K_i(V)^{\c}\to K_i(U\cap V)^{\c}\to\cdots.$$

$X\times\AA_R^1$ is again a toric variety. Moreover,
$U\times\AA^1_R$ and $V\times\AA^1_R$ are open toric subvarieties
covering $X\times\AA_R^1$ so that the underlying fans are products
of the corresponding fans by the same ray. Writing up the
corresponding Mayer-Vietoris sequences one obtains
$K_*(X)^{\c}=K_*(X\times\AA^1_R)^{\c}$, provided the same equality
holds for $U$ and $V$. Iterating the process we will produce toric
varieties whose fans contain less and less top dimensional cones.
But for single cones we are already done.

Conversely, using the same homotopy trick as in Step 2 in the
proof Theorem \ref{lgp} one can show that the equalities
$K_i(X\times\AA^1_R)^{\c}=K_i(X)^{\c}$ for affine toric varieties
$X$, whose cones are pointed, are in fact equivalent to Conjecture
\ref{conj}.
\end{proof}

\section{pyramidal extensions}\label{PYRAM}

\begin{definition}\label{pyr}
A monoid extension $M\subset N$ is called {\it pyramidal} if $M$
and $N$ are finitely generated normal monoids without nontrivial
units, $\gp(M)=\gp(N)$, and there is a representation of the form
$\Phi(N)=\Phi(M)\cup\delta$ for some rational pyramid $\delta$ so
that $\dim\Phi(N)=\dim\Phi(M)=\dim\delta$ and the polytope
$\Phi(M)$ meets $\delta$ in its base.
\end{definition}

Recall that a polytope is called {\it pyramid} if it is a convex
hull of its facet and a vertex not in this facet. The equivalent
reformulation of Conjecture \ref{conj} given in the lemma below
will be called {\it pyramidal descent}.

\begin{lemma}\label{piram}
Conjecture \ref{conj} is equivalent to the claim that for every
pyramidal extension $M\subset N$, a local regular ring $R$, an
element $x\in K_i(R[N_*])$ ($i\in\NN$) and a sequence
$\c=(c_1,\ldots)$ the inclusion $(c_1\cdots
c_j)_*(x)\in\I{\big(}K_i(R[M_*])\to K_i(R[N_*]){\big)}$ holds for
$j\gg0$.
\end{lemma}

\begin{proof}
That the conjecture implies the claim is clear. Next observe that
the mentioned inclusions are equivalent to the surjectivity of the
homomorphisms $K_i(R[M^{\c}_*])\to K_i(R[N^{\c}_*])$ for all
pyramidal extensions $M\subset N$ and all sequences $\c$.

Let $N$ be a finitely generated normal monoid with trivial $U(N)$.
Consider any sequence $N=N_1,N_2,$ $N_3,\ldots$ of normal
submonoids of $N$, satisfying the condition -- for any $k\in\NN$
either $N_{k+1}\subset N_k$ is a pyramidal extension or
$N_k\subset N_{k+1}$. (Observe that $\gp(N_k)=\gp(N)$, $k\in\NN$.)
Then our condition yields surjectivity of the homomorphisms
$K_i(R[(N_k^{\c})_*])\to K_i(R[N^{\c}_*])$. Such a sequence of
monoids will be called {\it admissible}.

According to Lemma 2.8 in \cite{Gu1} for any rational interior
point $\xi\in\Phi(N)$ and its any neighborhood $\xi\in
U\subset\Phi(N)$ there is an admissible sequence of monoids
$N,N_2,N_3,\ldots$ such that $\Phi(N_k)\subset U$ for all
sufficiently big indices $k$.

By Lemma \ref{int} there exists a rational simplex
$\Delta\subset\int(\Phi(N))$ such that $N(\Delta)$ is a free
submonoid of $N_*$. These observations altogether show that
(assuming the pyramidal descent) the homomorphism
$K_i(R[N(\Delta)^{\c}_*])\to K_i(R[N^{\c}_*])$ is surjective. But
then $K_i(R[N(\Delta)^{\c}])\to K_i(R[N^{\c}_*])$ is a surjection
as well. Since $N(\Delta)^{\c}$ is a filtered union of free
monoids Quillen's theorem implies $K_i(R[N^{\c}_*])=K_i(R)$. By
Lemma \ref{patch} we are done.
\end{proof}

\section{Approximations by polarized monoids}\label{POLAR}

The results of this section are refined versions of the results
from \S1.3 in \cite{Gu2}.

\begin{definition}\label{polar}
A {\it polarized monoid} $N$ is a triple $(t,\Gamma,N)$ where $N$
is a finitely generated normal monoid with trivial $U(N)$,
$\Gamma\subset\Phi(N)$ is a rational polytope and $t\in N$ is a
non-zero element -- the {\it pole}, such that the following hold:
\begin{itemize}
\item[(1)]
$C(N)=\RR_+t+\RR_+\Gamma$,
\item[(2)]
$N\cap(\RR_+t+\RR_+F)\approx\ZZ_+\times N(F)$ for any facet
$F\subset\Gamma$.
\end{itemize}
\end{definition}
In particular, in a polarized monoid $(t,\Gamma,N)$ one has that
$\RR_+t$ is an edge of $C(N)$, $t$ is the generator of
$N\cap\RR_+t\approx\ZZ_+$, and $t\notin\RR F$ (the hyperplane
spanned by $F$) for any facet $F\subset\Gamma$.

In a polarized monoid $(t,\Gamma,N)$ a facet $F\subset\Gamma$ is
called {\it positive} if $t$ and $\int(C(N)$ lie on the same side
relative to the hyperplane $\RR F\subset\RR\otimes\gp(M)$,
otherwise $F$ is called {\it negative}.

In what follows we will use the notation
$N^-=\ZZ_+(-t)+N(\Gamma)\subset\gp(N)$. Clearly, for a polarized
monoid $(t,\Gamma,N)$ the triple $(-t,\Gamma,N^-)$ is also a
polarized monoid. It will be called the {\it antipode} of
$(t,\Gamma,N)$.

The main approximation result for polarized monoids is as follows.

\begin{theorem}\label{approxB}
Let $M\subset N$ be a pyramidal extension of monoids. Then for any
natural number $s$ there exist systems of polarized monoids
$(t_\alpha,\Gamma_{1\alpha},N_{1\alpha})$,
$(t_\alpha,\Gamma_{2\alpha},N_{2\alpha})$, $\ldots$,
$(t_\alpha,\Gamma_{s\alpha},N_{s\alpha})$, $\alpha\in\NN$ such
that the following hold:
\begin{itemize}
\item[(a)]
$N_{1\alpha}\subset N_{2\alpha}\subset\cdots\subset N_{s\alpha}
\subset N^{\c}_*$ and
$\Gamma_{1\alpha}\subset\Gamma_{2\alpha}\subset\cdots\subset
\Gamma_{s\alpha}\subset\int(\Phi(M))$ for any $\alpha$,
\item[(b)]
$\pm t_\alpha+(N_{l\alpha}(\Gamma_{l\alpha})\setminus\{0\})\subset
\int(N_{l+1\alpha}(\Gamma_{l+1\alpha}))$, for any $l\in[1,s-1]$
and any $\alpha$,
\item[(c)]
$M^{\c}_*$ is a filtered union of the
$N_{1\alpha}(\Gamma_{1\alpha})$ and $N^{\c}_*$ is a filtered union
of the $N_{1\alpha}$.
\end{itemize}
\end{theorem}

The condition (c) automatically implies that $M^{\c}_*$ is a
filtered union of any of the  systems
$\{N_{2\alpha}(\Gamma_{2\alpha})\}_\alpha$,$\ldots$,$\{N_{s\alpha}
(\Gamma_{s\alpha})\}_\alpha$ and $N^{\c}_*$ is a filtered union of
any of the systems
$\{N_{2\alpha}\}_\alpha$,$\ldots$,$\{N_{s\alpha}\}_\alpha$. We
also remark that by dropping the condition (b) and assuming
$c_1=c_2=\cdots$ one exactly gets {\it Approximation Theorem B}
from \cite{Gu2}.

In the proof we will use

\begin{lemma}\label{approxA}
Let $L$ be a monoid such that $\Phi(L)$ is a simplex. Then
$L^{\c}_*$ is a filtered union of free monoids
\end{lemma}

The monoids of this type are usually called {\it simplicial} in
the toric literature. Notice, that monoids of the type $L^{\c}$
for $L$ simplicial are in general {\it not} filtered unions of
free monoids, an example --
${\big(}\ZZ_+(2,0)+\ZZ_+(1,1)+\ZZ_+(0,2){\big)}^{\bar 3}$ where
$\bar 3=(3,3,\ldots)$.

Lemma \ref{approxA} is exactly Approximation Theorem A from
\cite{Gu2} in the special case when $\c$ is a constant sequence.
The proof of the general claim makes no difference.

Fisrt let us show how Lemma \ref{approxA} proves Conjecture
\ref{conj} in the simplicial case.

\begin{theorem}\label{simplok}
Conjecture \ref{conj} is true for arbitrary simplicial monoid $M$.
\end{theorem}

\begin{proof} Assume $R$ is a regular ring. By Lemma \ref{exc}
$R\int(M^{\c})$ satisfies excision. On the other hand the exact
sequence $$ 0\to R\int(M^{\c})\to R[M^{\c}_*]\to
R[M^{\c}_*]/R\int(M^{\c})\to0 $$ and Lemma \ref{approxA} show that
$K_i(R\int(M^{\c}))=0$ for $i\in\ZZ$. (Here we use the
$K$-homotopy invariance of $R$.) Since all faces of $\Phi(M)$ are
also simplices the same arguments as in the proof of Lemma
\ref{inter} provide the desired process of ``annihilating'' the
interior monoids.
\end{proof}

{\it Proof of Theorem \ref{approxB}}. Fix an index $j\in\NN$ and
finite subsets $W\subset\int(\Phi(N))$, $W'\subset\int(\Phi(M))$.
We need to show that there are polarized monoids
$(t,\Gamma_1,N_1)$,$\ldots$, $(t,\Gamma_s,N_s)$ for which the
conditions (a) and (b) are satisfied and, moreover,
$W\subset\Phi(N_1)$, $W'\subset\Gamma_1$ and $$
G_j:=\frac{\gp(N)}{c_1\cdots c_j}\subset\gp(N_1). $$ Let $v$ be
the vertex of $\Phi(N)$ not in $\Phi(M)$. We can choose a rational
point $v'\in\int(\Phi(N))\setminus\Phi(M)$ close to $v$ so that it
does not belong to any of the hyperplanes $\RR F$
$(\RR\otimes\gp(N))$, where $F$ runs through the facets of
$\Phi(M)$.

Consider any rational simplex $\Delta$ in the affine hull of
$\Phi(N)$ which has one vertex at $v'$ and contains $W$ in its
interior (such exists). We have the normal monoid
$$
L=\{n\in\gp(N),\ n\not=0|\ \
\RR_+n\cap\Delta\not=\emptyset\}\cup\{0\}.
$$
By Lemma \ref{Gord} $L$ is a {\it finitely generated} normal
monoid and by Lemma \ref{int}(c) $\gp(N)=\gp(L)$. By Lemma
\ref{approxA} there exists a free submonoid $L_0\subset L^{\c}_*$
such that $G_j\subset\gp(L_0)$, $W\subset\int(\Phi(L_0))$ and the
simplex $\Phi(L_0)$ has a vertex $v_0$ close to $v'$, not in the
affine span of any of the facets $F\subset\Phi(M)$.

Put $\Delta_0=\Phi(L_0)$ and let $\delta_0\subset\Delta_0$ be the
facet opposite to $v_0$. The free generator of
$L_0(\{v_0\})\approx\ZZ_+$ will be denoted by $\tau$.

Consider the polar projection $\pi:\Phi(M)\to\Aff(\delta_0)$ into
the affine hull of $\delta_0$ with respect to the pole $v_0$. Then
$\dim\pi(F)=\dim F$ ($=\dim\delta_0$) for any facet
$F\subset\Phi(M)$. By Lemma \ref{int}(c) for any facet
$F\subset\Phi(M)$ there is a free basis
$$
l_{\pi(F),1},\ldots,l_{\pi(F),\rank N-1}\subset L(\pi(F))
$$
of the group $\gp{\big(}L(\delta_0){\big)}\approx\ZZ^{\rank N-1}$.
Fix such bases arbitrarily. We can also fix a natural number
$j'>j$ such that
$$
\tau_k=\frac\tau{c_{j'}c_{j'+1}\cdots c_{j'+k}}\in N^{\c}_*,\ \
k\in\NN.
$$
For a natural number $k$ and a facet $F\subset\Phi(M)$ we have the
systems
$$
\{\tau_k,l_{\pi(F),1}+a_1\tau_k,\ldots,l_{\pi(F), \rank
N-1}+a_{\rank N-1}\tau_k\}\subset N^{\c}_*
$$
where $a_1,\ldots,a_{\rank N-1}\in\ZZ_+$. The free monoids
$N(F,k,a_1,\ldots,a_{\rank N-1})$ they generate all have the same
groups of differences, provided $k$ is fixed. Denote by $G_k$ this
common group of differences.

We let $\delta(F,k,a_1,\ldots,a_{\rank N-1})\subset
\Phi(N(F,k,a_1,\ldots,a_{\rank N-1}))$ denote the facet opposite
to $v_0$.

It is an elementary geometric observation that if $k$ is
sufficiently big then the rational points
$\RR_+(l_{\pi(F),i}+a_i\tau_k)\cap\Phi(N)$, $i\in[1,\rank N-1]$
move ``sufficiently slowly'' towards the point $v_0$ when the
$a_i$ run through $\ZZ_+$. Moreover, all the time these rational
points remain correspondingly in the segments
$[\RR_+l_{\pi(F),i}\cap\Phi(N),\RR_+\tau_k\cap\Phi(N)]$,
$i\in[1,\rank N-1]$. In particular, we can choose $k$ big enough
and $a_{F,1},\ldots,a_{F,\rank L-1}\in\ZZ_+$ in such a way that
the affine hulls of the simplices
$\delta(F,k,a_{F,1},\ldots,a_{F,\rank L-1})$ bound a (convex)
subpolytope $\Gamma_1\subset\int(\Phi(M))$ with the same number of
facets as $\Phi(M)$ such that $W'\subset\Gamma_1$, while the
convex hull of $\Gamma_1$ and $v_0$ contains $W$. Moreover, by
approximating the facets of $\Phi(M)$ with a sufficient precision,
we can also achieve that $v_0$ is not in the affine hull of any of
the facets of $\Gamma_1$.

It is an easy exercise to show that the triple
$(\tau_k,\Gamma_1,N_1)$ is a polarized monoid, where
$N_1=\RR_+\conv(\Gamma_1,v_0)\cap G_k$. For instance, the
normality of $N_1$ follows from the facts that $\Phi(N_k)$ is a
union of the pyramids with vertex at $v_0$ and bases the positive
facets of $\Gamma_1$, that these pyramids define normal submonoids
of $N_1$ (they are free extensions of the normal monoids
corresponding to these facets), and that their groups of
differences are all the same.

Moreover, $\Gamma_1$ can even be chosen in such a way that there
are $2s-1$ polytopes $\Gamma_2,\ldots,\Gamma_{2s}$, obtained in
the same way as $\Gamma_1$, for which the following hold
\begin{itemize}
\item
$\Gamma_1\subset\int(\Gamma_2)$, $\ldots$, $\Gamma_{2s-1}\subset
\int(\Gamma_{2s})$ and $\Gamma_{2s}\subset\int(\Phi(M))$,
\item
$\Gamma_2,\ldots,\Gamma_{2s}$ possess all the aforementioned
properties of $\Gamma_1$.
\end{itemize}
In fact, one first ``refines'' the triple $(\tau_k,\Gamma_1,N_1)$
by passing to $\tau_{k'}$ with $k'>k$, leaving the polytope
$\Gamma_1$ untouched. The monoid $N_1$ is then changed by the
monoid $N'_1$ which is generated by $\tau_{k'}$ and the elements
of $N_1$ living on the positive facets of $\Gamma_1$. The point is
that this new $N'_1$ is again a polarized monoid and the addition
of $\pm\tau_{k'}$ to the elements of
$N'_1(\Gamma_1)\setminus\{0\}$ has sufficiently small effect on
the radial directions (provided $k'$ is big enough). In
particular, we can ``blow up'' the polytope $\Gamma_1$ suitably in
any prescribed number of steps, remaining inside $\int(\Phi(M))$.

It is immediate from Definition \ref{polar} that the $s$-tuple of
polarized monoids $$
(\tau_k,\Gamma_1,N_1),(\tau_k,\Gamma_3,N_3),\ldots,
(\tau_k,\Gamma_{2s-1},N_{2s-1}) $$ satisfies the desired
conditions.

Actually, we have even achieved
$\gp(N_1)=\gp(N_3)=\cdots=\gp(N_{2s-1})$. \qed

\section{The scheme $X$}\label{VARI}

From now on we will mainly work over a coefficient field $k$, not
necessarily algebraically closed.

Let $(t,\Gamma,N)$ denote a polarized monoid and
$(t\1,\Gamma,N^-)$ be its antipode (in the multiplicative
notation). The non-affine scheme $X$ is defined by gluing
$\Spec(k[N])$ and $\Spec(k[N^-])$ along their open subscheme
$\Spec(k[t\1N])$. (Here the equality $t\1N=tN^-$ is used.)

\begin{remark}\label{toric}
(a) In other words $X$ is a toric variety over $k$ whose fan
consists of 2 maximal equidimensional cones, sharing a facet.
Namely, we mean the cones $C(N)^*,C(N^-)^*\subset(\RR^{\rank
N})^*$, where $-^*$ refers to ``dual''. We also have the equality
$C(N)^*\cup C(N^-)^*=C(N(\Gamma))^*$. In particular, $X$ is
quasiprojective over $k$. In fact, the fan of $X$ admits a
completion to a {\it projective} fan (the one defining a
projective toric variety) -- an elementary geometric observation.
Therefore, $X$ is an open subscheme of a projective toric variety
(see \cite[Ch. 2]{O}).

(b) Formal similarities of $X$ with a projective line allow us to
use in \S\ref{PDESCENT} some of Quillen's ideas from his
computation of the $K$-theory of a projective line
\cite[\S8.3]{Q1}. However, it is not until \S\ref{PYDESC} that we
use peculiar properties of polarized monoids.
\end{remark}

The category of sheaves of locally free $\O_X$-modules will be
denoted by $\PP(X)$. Since $X$ is connected any object of $\PP(X)$
has a constant rank.

The objects of $\PP(X)$ can be thought of as triples
$(P,P^-,\Theta)$ where $P\in\Pp(k[N])$, $P^-\in\Pp(k[N^-])$, and
$\Theta:t\1P\to tP^-$ is an isomorphism of $k[t\1N]$-modules. In
this terminology a morphism
$(P_1,P_1^-,\Theta_1)\to(P_2,P_2^-,\Theta_2)$ is just a pair of
morphisms $f^+:P_1\to P_2$ and $f^-:P_1^-\to P_2^-$
correspondingly in the categories $\Pp(k[N])$ and $\Pp(k[N^-])$,
such that $f^-\circ\Theta_1=\Theta_2\circ f$.

Given an object $(P,P^-,\Theta)\in\PP(X)$. By \cite{Gu1} $P$
($P^-$) is free $k[N]$-module ($k[N^-]$-module). By fixing bases
in $P$ and $P^-$ we can associate to $\Theta$ a matrix $\theta\in
GL_{\rank P}(k[t\1N])$. We will say that $\Theta$ is represented
by $\theta$.

The following is obvious.

\begin{lemma}\label{matrix}
Assume $(P_1,P_1^-,\Theta_1),(P_2,P_2^-,\Theta_2)\in\PP(X)$.
\begin{itemize}
\item[(a)]
If $\Theta_1$ is represented by a matrix $\theta_1$ then the set
of matrices representing $\Theta_1$ is exactly
$\{\tau\theta_1\sigma\ |\ \tau\in GL_{\rank P}(k[N]),\ \sigma\in
GL_{\rank P}(k[N^-])\}$.
\item[(b)]
$(P_1,P_1^-,\Theta_1)\approx(P_2,P_2^-,\Theta_2)$ if and only if
$\rank P_1=\rank P_2$ and for some (equivalently, for all)
matrices $\theta_1$ and $\theta_2$, representing correspondingly
$\Theta_1$ and $\Theta_2$, there exist $\tau$ and $\sigma$ as
above such that $\tau\theta_1\sigma=\theta_2$.
\end{itemize}
\end{lemma}

We have the natural augmentation $k[t\1N]\to k[t^{\pm1}]$, defined
by $t\1N\setminus U(t^{-1}N)\to0$. For a matrix $\beta$, defined
over $k[t\1N]$, its image under this augmentation will be denoted
by $\beta(0)$.

\begin{definition}\label{PPR0}
$\PP(X)^0\subset\PP(X)$ is the full subcategory of the objects
$(P,P^-,\Theta)$ such that $\Theta$ is represented by a matrix
$\theta$ satisfying the condition $[\theta(0)\1\cdot\theta]=0\in
K_1(k[t\1N])$.
\end{definition}

Any natural number $c$ gives rise to an endomorphism of $X$ which
on the affine charts
$\Spec(k[N]),\Spec(k[N^-])\subset\Spec(k[t\1N])$ is given by
$k[-^c]$. The induced endomorphism of $K_i(X)$ ($i\in\NN$) will be
denoted by $c_*$.

\begin{lemma}\label{element}
$\PP(X)^0$ is closed under extensions in $\PP(X)$ and for an
element $z\in K_i(X)$ we have $(c_1\cdots
c_j)_*(z)\in\I{\big(}K_i(\PP(X)^0)\to K_i(X){\big)}$ provided
$j\gg0$.
\end{lemma}

In the proof we will use the main result of \cite{Gu2} (the stable
version):

\begin{theorem}\label{mainG2}
$SK_1(k[L^{\c}])=0$ for any monoid $L$.
\end{theorem}

The stronger result $SL_n(k[L^{\c}])=E_n(k[L^{\c}])$ for $n\geq3$
is proved in \cite[Theorem 2.1]{Gu2} when $c_1=c_2=\cdots$. But
the same claim remains true for arbitrary $\c$. In fact, the only
place in \cite{Gu2} that uses the sequence $(c,c,\ldots)$ is
Approximation Theorem B there. But Theorem \ref{approxB} is a
refined version for arbitrary $\c$.

{\it Proof of Lemma \ref{element}.} The closedness under
extensions of $\PP(X)^0$ follows from the following observations.
Assume $(P_1,P_1^-,\Theta_1),(P_2,P_2^-,\Theta_2)\in\PP(X)^0$ and
$\theta_1\in GL_{\rank P_1}(k[t\1N])$ and $\theta_2\in GL_{\rank
P_2}(k[t\1N])$ represent correspondingly $\Theta_1$ and
$\Theta_2$. Then for an exact sequence in $\PP(X)$ of the type
$$
0\to(P_1,P_1^-,\Theta_1)\to(P_3,P_3^-,\Theta_3)\to(P_2,P_2^-,\Theta_2)\to0
$$
there is a representation of $\Theta_3$ by a matrix of the form
$$
\theta_3=
\begin{pmatrix}
\theta_1&\theta'\\
0&\theta_2
\end{pmatrix}\in GL_{\rank P_1+\rank P_2}(k[t\1N])
$$
where $\theta'$ is a $(\rank P_1)\times(\rank P_2)$-matrix over
$k[t\1N]$. It is clear that $[\theta_3(0)\1\theta_3]=0\in
K_1(k[t\1N])$ provided $[\theta_1(0)\1\theta_1]=0$ and
$[\theta_2(0)\1\theta_2]=0$ in $K_1(k[t\1N])$.

Let a continuous map $\zeta:S^{i+1}\to BQ\PP(X)$ represent $z$.
($S^{i+1}$ is $(i+1)$-sphere.) There are finitely many simplices
in the simplicial complex $BQ\PP(X)$ which intersect $\I(\zeta)$.
Let $(P_\lambda,P_\lambda^-,\Theta_\lambda)$, $\lambda\in\Lambda$
be the vertices of these simplices and the matrices
$\theta_\lambda\in GL_{\rank P_\lambda}(k([t\1N]))$ represent the
$\Theta_\lambda$. Since $U(k[t\1N])=U(k[t^{\pm1}])$, Theorem
\ref{mainG2} implies $[(c_1\cdots
c_j)_\#{\big(}\theta_\lambda(0)\1\cdot \theta_\lambda{\big)}]=0$
in $K_1(k[t\1N])$, $\lambda\in\Lambda$, provided $j\gg0$. (Here
$c_{\#}$ refers to the corresponding endomorphism of
$GL(k[t\1N])$.) \qed

\section{Polarization descent}\label{PDESCENT}

In this sections we use the same notation as in \S\ref{VARI}.

For a natural number $r$ we denote by $\PP(X)(r)^0\subset\PP(X)^0$
the full subcategory of the objects whose ranks are nonnegative
multiples of $r$. It is clear that $\PP(X)(r)^0$ is closed under
extensions in $\PP(X)^0$.

\begin{definition}\label{diff}
Assume $r\in\NN$ and $a,b\in\ZZ\cup\{-\infty,+\infty\}$, $a\leq
b$. Then
$$
\PP(X)(r)_{a,b}^0\subset\PP(X)(r)^0
$$
is the full subcategory of the objects $(P,P^-,\Theta)$ such that
$\Theta$ can be represented by a matrix $\theta\in GL_{\rank
P}(k[t\1N])$ satisfying the conditions:
\begin{itemize}
\item[(1)]
$$
\theta(0)=
\begin{pmatrix}
t^{u_1}&0&\ldots&0\\
0&t^{u_2}&\ldots&0\\
.&.&.&.\\
0&0&\ldots&t^{u_{\rank P}}
\end{pmatrix}
$$
for some integral numbers $a\leq u_1,\ldots,u_{\rank P}\leq b$,
\item[(2)]
$[\theta(0)\1\theta]=0\in K_1(k[t\1N])$.
\end{itemize}
\end{definition}

The {\it polarization} of an object of
$\PP(X)(r)^0_{-\infty,+\infty}$ is by definition the smallest
segment $[a,b]$ such that this object is in $\PP(X)(r)^0_{a,b}$.
The polarization is well defined because
${\big(}\O_{\Pp^1}(u_1)\oplus\cdots\oplus\O_{\Pp^1}(u_n)\approx
\O_{\Pp^1}(v_1)\oplus\cdots\oplus\O_{\Pp^1}(v_m){\big)}$
$\Rightarrow$ ${\big(}n=m\ \text{and the}\ u\ \text{coincide with
the}\ v\ \text{up to permutation}{\big)}$. Here $\Pp^1=\Pp^1_k$
denotes the projective line over $k$ and
$\O_{\Pp^1}(u)=(k[t],k[t^{-1}],\ \text{multiplication by}\
t^{-u})$.

\begin{lemma}\label{poldesc}
$\PP(X)(r)^0_{0,1}$ is closed under extensions in $\PP(X)(r)^0$
and the embedding $\PP(X)(r)^0_{0,1}\subset\PP(X)(r)^0$ induces
isomorphisms on $K$-groups.
\end{lemma}

\begin{proof}
We will only treat the case $r=1$. The general case makes no
difference.

{\it Step 1.} Here we prove that $\PP(X)^0_{-\infty,+\infty}=
\PP(X)^0$.

Consider an object $(P,P-,\Theta)\in\PP(X)^0$. By Definition
\ref{PPR0} $\Theta$ can be represented by a matrix $\theta\in
GL_{\rank P}(k[t\1N])$ such that $[\theta(0)\1\theta]=0\in
K_1(k[t\1N])$. By Grothendieck's theorem \cite{Gro} on vector
bundles on projective lines (proved by elementary methods for
arbitrary fields in \cite{HMa}) there are $\sigma\in GL_{\rank
P}(k[t])$ and $\tau\in GL_{\rank P}(k[t\1])$ such that
$$ \sigma\theta(0)\tau=(\sigma\theta\tau)(0)
\begin{pmatrix}
t^{u_1}&0&\ldots&0\\
0&t^{u_2}&\ldots&0\\
.&.&.&.\\
0&0&\ldots&t^{u_{\rank P}}
\end{pmatrix}
$$ for some $u_1,\ldots,u_{\rank P}\in\ZZ$. We have
$$((\sigma\theta\tau)(0))^{-1}(\sigma\theta\tau)
=\tau\1\theta(0)\1\theta\tau=0\in K_1(k[t\1N]). $$ By Lemma
\ref{matrix} we are done.

{\it Step 2.} We prove that $\PP(X)^0_{a,b}$ is closed under
extensions in $\PP(X)^0$ for arbitrary
$a,b\in\ZZ\cup\{\pm\infty\}$, $a\leq b$.

By the previous step the reduction modulo $k(t^{-1}N\setminus
t^{\ZZ})\subset k[t\1N]$ shows that it is enough to prove the
following

\

{\bf Claim.} Assume we are given a commutative diagram of free
$k[t^{\pm1}]$-modules
$$
\xymatrix{ \zto k[t^{\pm1}]^m\ar[d]^\alpha\ar[r]^-{f^+}
&k[t^{\pm1}]^{m+n}        \ar[d]^\beta \ar[r]^-{g^+}
&k[t^{\pm1}]^n\toz\ar[d]^\gamma\\
\zto k[t^{\pm1}]^m\ar[r]_-{f^-} &k[t^{\pm1}]^{m+n}\ar[r]_-{g^-}
&k[t^{\pm1}]^n\makebox[0pt]\toz }
$$
with exact rows, where $f^+$ and $g^+$ are defined over $k[t]$ and
$f^-$ and $g^-$ are defined over $k[t\1]$. Assume, further,
$\alpha$, $\beta$ and $\gamma$ are $k[t^{\pm1}]$ isomorphisms
determined correspondingly by matrices of the type
$$
\begin{pmatrix}
t^{u_1}&0&\ldots&0\\
0&t^{u_2}&\ldots&0\\
.&.&.&.\\
0&0&\ldots&t^{u_m}
\end{pmatrix},\ \
\begin{pmatrix}
t^{v_1}&0&\ldots&0\\
0&t^{v_2}&\ldots&0\\
.&.&.&.\\
0&0&\ldots&t^{v_{m+n}}
\end{pmatrix}\ \text{and}\
\begin{pmatrix}
t^{w_1}&0&\ldots&0\\
0&t^{w_2}&\ldots&0\\
.&.&.&.\\
0&0&\ldots&t^{w_n}
\end{pmatrix}.
$$
Then $(u_1,\ldots,u_m,w_1,\ldots,w_n\in[a,b])\Rightarrow
(v_1,\ldots,v_{m+n}\in[a,b])$.

\

We introduce the following notation:

(i) the matrices of $f^+$ and $f^-$ (w.r.t. the standard bases)
are correspondingly $(f_{ij}^+)$, $i\in[1,m]$, $j\in[1,m+n]$,
$f_{ij}^+\in k[t]$ and $(f_{ij}^-)$, $i\in[1,m]$, $j\in[1,m+n]$,
$f_{ij}^-\in k[t^-]$,

(ii) the matrices of $g^+$ and $g^-$  are correspondingly
$(g_{jl}^+)$, $j\in[1,m+n]$, $l\in[1,n]$, $g_{jl}^+\in k[t]$ and
$(g_{jl}^-)$, $j\in[1,m+n]$, $l\in[1,n]$, $g_{jl}^-\in k[t\1]$.

The commutative diagram above implies $$
f_{ij}^+t^{v_j}=f_{ij}^-t^{u_i}\ \text{for all}\
i,j\eqno(\ref{poldesc})_1 $$ and $$
g_{jl}^+t^{w_l}=g_{jl}^-t^{v_j}\ \text{for all}\
j,l.\eqno(\ref{poldesc})_2 $$

Fix any index $j\in[1,m+n]$. If there exists $i\in[1,m]$ such that
$f_{ij}^+\not=0$ then (\ref{poldesc})$_1$ implies $u_i\geq v_j$.
Now consider the case when $f_{ij}^+=f_{ij}^-=0$ for all
$i\in[1,m]$. Then our diagram splits as follows $$ \xymatrix{ \zto
k[t^{\pm1}]^m \ar[d]^\alpha\ar[r]^-{f^+}
&k[t^{\pm1}]^{m+n-1}\oplus k[t^{\pm1}] \ar[d]^\beta
\ar[r]^-{g_1^+\oplus g_2^+} &k[t^{\pm1}]^{n-1}\oplus
k[t^{\pm1}]\toz \ar[d]^\gamma\\ \zto k[t^{\pm1}]^m\ar[r]_-{f^-}
&k[t^{\pm1}]^{m+n-1}\oplus k[t^{\pm1}]\ar[r]_-{g_1^-\oplus g_2^-}
&k[t^{\pm1}]^{n-1}\oplus k[t^{\pm1}]\toz } $$ where $$ 0\to
k[t^{\pm1}]^m\xto{f^{\pm}}k[t^{\pm1}]^{m+n-1}\xto{g_1^{\pm}}
k[t^{\pm1}]^{n-1}\to0$$ are the exact sequences of free
$k[t^{\pm1}]$-modules, obtained by deleting the $j$-th direct
summand in $k[t^{\pm1}]^{m+n}$, and the mappings
$g_2^+,g_2^-:R[t^{\pm1}]\to k[t^{\pm1}]$ establish an isomorphism
between the $\O_{\Pp^1}$-sheaves $\O(-v_j)$ and $\O(-w_l)$ for
some $l\in[1,n]$. But in this case $v_j=w_l$. Therefore, in all
cases we have $v_j\leq\max_{i,l}(u_i,w_l)$.

Using (\ref{poldesc})$_2$ and dual arguments we get
$v_j\geq\min_{i,l}(u_i,w_l)$.

{\it Step 3.} By Step 1 we have the filtered union representation
$$
\PP(X)^0=\bigcup_{a=0}^{-\infty}\PP(X)^0_{a,+\infty}.
$$
Therefore, if we knew that the embeddings
$\PP(X)^0_{a,+\infty}\subset\PP(X)^0_{a-1,+\infty}$, $a\leq0$,
induce isomorphisms on $K$-groups we could conclude that the
embedding $\PP(X)^0_{0,+\infty}\subset\PP(X)^0$ also induces
isomorphisms on $K$-groups.

Fix a nonpositive integral number $a$ and consider the two exact
functors
$$
F_1,F_2:\PP(X)^0_{a-1,+\infty}\to\PP(X)^0_{a,+\infty}
$$
given correspondingly by
$$
F_1((P,P^-,\Theta))=(P,P^-,t\Theta)\oplus(P,P^-,t\Theta),\ \
F_1((f^+,f^-))=(f^+,f^-)
$$
and
$$
F_2((P,P^-,\Theta))=(P,P^-,t^2\Theta),\ \
F_1((f^+,f^-))=(f^+,f^-).
$$
The exact ``Koszul sequence'' (where we use matrix theoretical
notation)
$$
\xymatrix@!C{ \zto k[t\1N]^n               \ar[d]_{t^2\Theta}
\ar[r]^-{(t,-1)^{\oplus n}} &k[t\1N]^n\oplus
k[t\1N]^n\ar[d]_{t\Theta\oplus t\Theta}\ar[r]^-{{\binom1t}^{\oplus
n}} &k[t\1N]^n\toz               \ar[d]_\Theta
\\
\zto k[t\1N]^n               \ar[r]_-{(1,-t\1)^{\oplus n}}
&k[t\1N]^n\oplus k[t\1N]^n\ar[r]_-{{\binom{t\1}1}^{\oplus n}}
&k[t\1N]^n\toz } $$ implies the exact sequence of exact functors
$0\to\iota\circ F_2\to\iota\circ F_1\to{\bf 1}\to0$, where
$\iota:\PP(X)^0_{a,+\infty}\to\PP(X)^0_{a-1,+\infty}$ is the
identity embedding and $\bf 1$ is the identity endofuntor of
$\PP_R(X)^0_{a-1,+\infty}$. By \cite[\S3]{Q1} we get
$\iota_*\circ{\big(}(F_1)_*-(F_2)_*{\big)}= {\bf 1}_*$ for the
corresponding $K$-group homomorphisms.

It remains to show that
${\big(}(F_1)_*-(F_2)_*{\big)}\circ\iota_*={\bf 1}_*$, the right
hand side denoting the identity endomorphism of the corresponding
$K$-group of $\PP(X)^0_{a,+\infty}$. But this equality is derived
by literally the same arguments once we observe the exact ``Koszul
sequence'' $0\to F_2\circ\iota\to F_1\circ\iota\to{\bf 1}\to0$ of
endofunctors  of $\PP(X)^0_{a,+\infty}$.

{\it Step 4.} Because of the filtered union representation
$$
\PP(X)^0_{0,+\infty}=\bigcup_{b=1}^{+\infty}\PP(X)^0_{0,b}
$$
the previous step shows that we only need that the embeddings
$\PP(X)^0_{0,b}\subset\PP(X)^0_{0,b+1}$, $b\geq1$, induce
isomorphisms on $K$-groups.

Assume $b\in[1,+\infty)$. Consider an object
$(P,P^-,\Theta)\in\PP(X)^0_{0,b+1}$ and fix a matrix $\theta$
satisfying the conditions in Definition \ref{diff} (with respect
to $a=0$ and $b=b+1$). We let $u_1,\ldots,u_n\in[0,b+1]$ be the
corresponding numbers, where $n=\rank P$.

Without loss of generality we can assume $u_1,\ldots,u_l=b+1$ and
$u_{l+1},\ldots,u_n\in[0,b]$ for some $l\in[0,n]$ (none of the
values $l=0$ and $l=n$ being excluded).

Consider the matrix
$$
\rho=
\begin{pmatrix}
t{\bf1}_l&0\\
0&{\bf1}_{n-l}
\end{pmatrix}.
$$
We have the following (non-functorial) commutative diagram with
exact rows
$$
\xymatrix{ \zto k[t\1N]^n\ar[d]^\theta\ar[r]^-{(\rho,-{\bf1})}
&k[t\1N]^n\oplus k[t\1N]^n
                 \ar[d]^{\spmatrix{\rho\1\theta&0\\0&\theta\rho\1}}
                              \ar[r]^-{\binom{\bf1}\rho}
&k[t\1N]^n\toz\ar[d]^{\rho\1{\theta}\rho\1}
\\
\zto k[t\1N]^n               \ar[r]_-{({\bf1},-\rho\1)}
&k[t\1N]^n\oplus k[t\1N]^n\ar[r]_-{\binom{\rho\1}{\bf1}}
&k[t\1N]^n\toz }
$$
where ${\bf1}={\bf1}_{n\times n}$ and the matrices at the vertical
rows refer to the corresponding maps.

Observe that the matrices at the second and third vertical rows
satisfy the conditions (1) and (2) in Definition \ref{diff} so
that we get objects of polarization $[0,b]$ (one uses that
$\rho(0)=\rho$). Therefore, any object
$(P,P^-,\Theta)\in\PP(X)^0_{0,b+1}$ admits a ``co-resolution''
$$
0\to(P,P^-,\Theta)\to(P_1,P^-_1,\Theta_1)\to(P_2,P^-_2,\Theta_2)\to0
$$
whose second and third terms belong to $\PP(X)^0_{0,b}$. In other
words, for the dual categories
$$
\PP^b={\big(}\PP(X)^0_{0,b}{\big)}^{\text{op}}\ \ \text{and}\ \
\PP^{b+1}={\big(}\PP(X)^0_{0,b+1}{\big)}^{\text{op}}
$$
any object ${\frak p}\in\PP^{b+1}$ admits a resolution $0\to{\frak
p}_2\to{\frak p}_1\to{\frak p}\to0$ where ${\frak p}_1,{\frak
p}_2\in\PP^b$. Since $K$-theories of dual exact categories are the
same it suffices to show that the embedding
$\PP^b\subset\PP^{b+1}$ induces isomorphisms on $K$-groups. In
view of what has been said above we see that Resolution Theorem
\cite[\S4]{Q1} applies once we show the implication: if
$0\to{\frak p}_2\to{\frak p}_1\to{\frak p}_0\to0$ is exact in
$\PP^{b+1}$ and ${\frak p}_1\in\PP^b$ then ${\frak p}_2\in\PP^b$.

Returning to the original categories, we claim that for any exact
sequence
$$
0\to(P_0,P^-_0,\Theta_0)\to(P_1,P^-_1,\Theta_1)\to(P_2,P^-_2,\Theta_2)\to0
$$
in $\PP(X)^0_{0,b+1}$ the following implication holds
$$
{\big(}(P_1,P^-_1,\Theta_1)\in\PP(X)^0_{0,b}{\big)}\Rightarrow
{\big(}(P_2,P^-_2,\Theta_2)\in\PP(X)^0_{0,b}{\big)}.
$$

Introducing the notation as in Claim in Step 2 this implication
rewrites as follows:
$$
(v_1,\ldots,v_{m+n}<b+1)\Rightarrow(w_1,\ldots,w_n<b+1).
$$
First observe that for any $h\in[1,n]$ there exists
$j_h\in[1,m+n]$ such that $g_{j_hh}^+\not=0$. In fact, assuming to
the contrary that such $j_h$ does not exist, the mapping $g^+$
(and, therefore, $g^-$ too) would not be surjective, which is
excluded. By (\ref{poldesc})$_2$ we have
$$
g_{j_hh}^+t^{w_h}=g_{j_hh}^-t^{v_{j_h}},\ \ h\in[1,n].
$$
In particular, $w_h\leq v_{j_h}$. This completes the proof.
\end{proof}

Let $\PP(\Pp^1_k)_{0,1}$ denote the full subcategory of
$\PP(\Pp^1_k)$ whose objects are $\O^a\oplus\O(-1)^b$,
$a,b\in\ZZ_+$. It is clear from Claim in Step 2 in the proof above
that $\PP(\Pp^1_k)_{0,1}$ is an exact category. For ${\mathcal
E},{\mathcal F}\in\PP(\Pp^1_k)_{0,1}$ the product ${\mathcal
F}^*\otimes{\mathcal E}$ is of type
$\O(-1)^a\oplus\O^b\oplus\O(1)^c$. Hence $0=H^1(\Pp^1_k,{\mathcal
F}^*\otimes{\mathcal E})=$ $\Ext({\mathcal F},{\mathcal E})$
(\cite[Ch.3 \S6]{Ha}), and we get

\begin{lemma}\label{split}
Any short exact sequence in $\PP(\Pp^1_k)_{0,1}$ splits.
\end{lemma}

\section{Almost pyramidal descent}\label{PYDESC}

The peculiar property of polarized monoids we need is Lemma
\ref{peculiar} below, yielding Claim A and Claim B in the proof of
Theorem \ref{pdescent}. We present it in the general form
involving arbitrary commutative rings.

First a word on notation. Let $(t,\Gamma,L)$ be a polarized
monoid. If $R$ is a local ring and $\mu\subset R$ is its maximal
ideal we let $\N\in\max R[L(\Gamma)]$ denote the maximal ideal
generated by $\mu$ and $L(\Gamma)\setminus\{1\}$. Also, for a
natural number $r$ we put
\begin{itemize}
\item
$G_r=E_r(R[t\1L]_{\N})\cap GL_r(R[t\1L]_{\N},\N R[t\1L]_{\N})$,
\item
$G^+_r=E_r(R[L]_{\N})\cap GL_r(R[L]_{\N},\N R[L]_{\N})$,
\item
$G^-_r=E_r(R[L^-]_{\N})\cap GL_r(R[L^-]_{\N},\N R[L^-]_{\N})$.
\end{itemize}
The subset $\{g^+g^-|\ g^+\in G_r^+,\ g^-\in G_r^-\}\subset G_r$
will be denoted by $G_r^+G^-_r$.

\begin{lemma}\label{peculiar}
$G_r=G_r^+G_r^-$ for any local ring $R$ and any natural number
$r\geq3$.
\end{lemma}

\begin{remark}\label{Sch}
Lemma \ref{peculiar} is a polarized version of the pivotal
technical fact in \cite{Su} (Proposition 5.6). It is proved in
\cite[Proposition 2.14, Step 4]{Gu2}. A generalization of Lemma
\ref{peculiar} to a relative case $I\subset R$ is obtained in
\cite{Sch}.
\end{remark}

Let $k$ be a field. Put
$$
\Lambda=\Lambda_{(t,\Gamma,L)}=\{(\phi_{uv})\}\subset
M_{2\times2}(k[L]),
$$
where
\begin{itemize}
\item
$\phi_{11},\phi_{22}\in k[L(\Gamma)]$,
\item
$\phi_{12}\in k[L(\Gamma)]\cap t\1k[L(\Gamma)]$,
\item
$\phi_{21}\in k[L(\Gamma)]+tk[L(\Gamma))]$.
\end{itemize}

Consider the subring $$\Lambda'=\Lambda'_{(t,\Gamma,L)}=
\{(\phi_{uv})|\ \phi_{21}(0)\in k\}\subset\Lambda.$$ (As usual,
$\phi(0)$ is the image of $\phi$ under the augmentation $k[L]\to
k[t]$.) Then $\Lambda$ differs from $\Lambda'$ in a single
monomial -- the pole $t$ (in the 21 position).

\begin{theorem}\label{pdescent}
Let $k$ be a field, $M\subset N$ be a pyramidal extension of
monoids and $z\in K_i(k[N_*])$ ($i\in\NN$). Assume Conjecture
\ref{conj} holds for all monoids of rank $<\rank N$  and all
coefficient fields. Then for every sequence $\c=(c_1,c_2,\ldots)$
there exists a polarized monoid $(t,\Gamma,L)$ such that:
\begin{itemize}
\item[(a)]
$L\subset N_*$ and $\Gamma\subset\int(\Phi(M))$,
\item[(b)]
$\Lambda'_{(t,\Gamma,L)}\subset M_{2\times2}(k[M_*])$,
\item[(c)]
$(c_1\cdots c_j)_*(z)\in \I{\big(}K_i(\Lambda_{(t,\Gamma,L)})\to
K_i(M_{2\times2}(k[N_*])) =K_i(k[N_*]){\big)}$ for $j\gg0$.
\end{itemize}
\end{theorem}

\begin{proof}
{\it Step 1.} Let $z(0)\in K_i(k)$ denote the image of $z$ under
the homomorphism $K_i(k[N_*])\to K_i(k)$, induced by
$N_*\setminus\{1\}\to0\in k$. Changing $z$ by $z-z(0)$ we can
without loss of generality assume $z(0)=0$.

Denote by $x$ the image of $z$ in $K_i(k[N^{\c}_*])$. Then
$x(0)=0$ w.r.t. the natural augmentation $k[N_*^{\c}]\to k$.

{\it Step 2.} By Theorem \ref{approxB} (and continuity of
algebraic $K$-functors) there are six polarized monoids
$(t,\Gamma_1,N_1)$,$\ldots$,$(t,\Gamma_6,N_6)$ such that
\begin{itemize}
\item[(i)]
$N_1\subset\cdots\subset N_6\subset N^{\c}_*$ and
$\Gamma_1\subset\cdots\subset\Gamma_6\subset\int(\Phi(M))$,
\item[(ii)]
$t^{\pm1}(N_l(\Gamma_l)\setminus\{1\})\subset\int(N_{l+1}(\Gamma_{l+1}))$,
for $l=[1,5]$.
\item[(iii)]
$x\in\I{\big(}K_i(k[N_1])\to K_i(k[N^{\c}_*]){\big)}$.
\end{itemize}
We will assume that $x$ is the image of $x_1\in K_i(k[N_1])$.

Next we show that one can also assume
\begin{itemize}
\item[(iv)]
$x_1\in\Ker{\big(}K_i(k[N_1])\to K_i(k[t\1N_1]){\big)}$.
\end{itemize}
Since $N_1,\ldots,N_6$ are finitely generated monoids there is an
index $j\in\NN$ such that $$ \frac{N_1}{c_j\cdots c_{j'}},\ldots,
\frac{N_6}{c_j\cdots c_{j'}}\subset N^{\c}_* $$ for any natural
number $j'\geq j$. Denote these six polarized monoids respectively
by $N_{1j'}$,  $\ldots$,$N_{6j'}$. They are naturally isomorphic
to $N_1,\ldots,N_6$. Clearly, the same conditions (i)--(iii) are
satisfied for them. In particular, if we show that
$x_1\in\Ker{\big(}K_i(k[N_1])\to K_i(k[(t\1N_1)^{\c'}]){\big)}$,
where $\c'=(c_j,c_{j+1},\ldots)$, then we achieve the validity of
all the conditions (i)--(iv) for the isomorphic system of monoids
$N_{1j'}$,$\ldots$,$N_{6j'}$ with $j'\gg0$.

By Corollary \ref{pyrapp} there is an intermediate monoid
$N_1\subset\ZZ_+t\times N_0\subset t\1N_1$, where $N_0\subset
t\1N_1$ is a normal submonoid for which $U(N_0)$ is trivial. In
view of the inclusions $k[N_1]\subset k[\ZZ_+t\times
N_0^{\c'}]=k[N_0^{\c'}][t]\subset k[(t\1N_1)^{\c'}]$ it suffices
to show that $x_1\in\Ker{\big(}K_i(k[N_1])\to
K_i(k[N_0^{\c'}][t]){\big)}$, By Lemma \ref{monic} this is
equivalent to the inclusion $x_1\in\Ker{\big(}K_i(k[N_1])\to
K_i(k(t)[N_0^{\c'}]){\big)}$. The hypothesis of the theorem
implies that the augmentation $k(t)[N_0^{\c'}]\to k(t)$,
$N_0^{\c'}\setminus\{1\}\to0$ induces the isomorphism
$K_i(k(t)[N_0^{\c'}])=K_i(k(t))$. Now the composite mapping
$k[N_1]\to k(t)[N_0^{\c'}]\to k(t)$ factors through $k[t]$.
Finally, $K_i(k[t])=K_i(k)$ and by Step 1 we are done.

{\it Step 3.} We let $X_1$ denote the non-affine toric variety
corresponding to $N_1$ in the sense of \S\ref{VARI}. By Theorem
\ref{MV} (and Remark \ref{quilwald}) we have the long exact
sequence $$ \cdots\to K_i(X_1)\to K_i(k[N_1])\oplus
K_i(k[N_1^-])\to K_i(k[t\1N_1])\to\cdots $$ Hence by (iv) in Step
2 there exists $z_1\in K_i(X_1)$ mapping to $x_1$.

For any natural number $r$ we have the commutative diagram of
exact categories and exact functors
$$
\xymatrix{
&&\PP(X_1(r))\ar[d]_\cap\\
\PP(X_1)              \ar[d]             \ar[r]_-\Delta\ar[urr]
&\PP(X_1)\times\cdots\times\PP(X_1)      \ar[r]_-\oplus
&\PP(X_1)             \ar[d]\\
\PP(k[N_1])                              \ar[r]_-\Delta
&\PP(k[N_1])\times\cdots\times\PP(k[N_1])\ar[r]_-\oplus
&\PP(k[N_1]) }
$$
where $\Delta$ refers to the diagonal embeddings and the number of
factors is $r$. Passing to $K$-groups we get the commutative
diagram
$$
\xymatrix{
&K_i(\PP(X_1)(r))\ar[d]\\
K_i(X_1)\ar[ur]\ar[r]_-{r\cdot}\ar[d]
&K_i(X_1)\ar[d]\\
K_i(k[N_1])\ar[r]_-{r\cdot} &K_i(k[N_1]) }.
$$

Therefore, for any natural number $r$ there is an element $z_r\in
K_i(\PP(X_1)(r))$ mapping to $rx_1$ under the homomorphism
$K_i(\PP(X_1)(r))\to K_i(k[N_1])$.

{\it Step 4.} As in Step 2 there exists an index $j\in\NN$ such
that $$ N_{1j'}=\frac{N_1}{c_j\cdots c_{j'}},\ldots,
N_{6j'}=\frac{N_6}{c_j\cdots c_{j'}}\subset
N^{\c}_*\eqno(\ref{pdescent})_1 $$ for any natural number $j'\geq
j$.

By Lemma \ref{element} for any natural number $r$ there exists
$j_r>j$ such that the element $(c_j\cdots c_{j_r})_*(z_r)$ is in
$\I{\big(}K_i(\PP(X_1)(r)^0)\to K_i(\PP(X_1)(r)){\big)}$. In
particular, by Step 3 we have
$rx\in\I{\big(K_i(\PP(X_{1j_r})(r)^0)}\to
K_i(k[N^{\c}_*]){\big)}$, where the mapping between the $K$-groups
is induced by the composite exact functor
$$
\PP(X_{1j_r})(r)^0\to\PP(X_{1j_r})(r)
\to\PP(k[N_{1j_r}])\to\PP(k[N^{\c}_*]).
$$
Here $X_{1j_r}$ is the scheme associated with $N_{1j_r}$. By Lemma
\ref{poldesc}
$$
rx_{\M}\in\I{\big(}K_i(\PP(X_{1j_r})(r)^0_{0,1})\to
K_i(k[N^{\c}_*]_{\M}){\big)},\eqno(\ref{pdescent})_2
$$
where the mapping between the $K$-groups is induced by the `upper
route' in the commutative diagram of exact categories and exact
functors
$$
\begin{aligned}
\xymatrix{ \PP(X_{1j_r})(r)\ar[r] &\PP(k[N_{1j_r}])\ar[r]
&\PP(k[N_*^{\c}])\ar[d]\\
\PP(X_{1j_r})(r)^0_{0,1}\ar[u]\ar[d] &
&\PP(k[N_*^{\c}]_{\M})\\
\PP((X_{1j_r})_{\N_{1j_r}})(r)\ar[r]
&\PP((X_{2j_r})_{\N_{2j_r}})(r)\ar[r]
&\PP((X_{4j_4})_{\N_{4j_r}})\ar[u] }
\end{aligned}
\eqno(\ref{pdescent})_3
$$
That the `lower route' in $(\ref{pdescent})_3$ is in fact possible
follows from $(\ref{pdescent})_1$ and the condition (i) in Step 2.
Here we have used the following notation:
$$
\M=k(M^{\c}_*\setminus\{1\})\in\max k[M^{\c}_*],
$$
$$
\N_{1j_r}=k{\big(}N_{1j_r}(\Gamma_1)\setminus\{1\}{\big)}\in\max
k[N_{1j_r}(\Gamma_1)],
$$
and the scheme $(X_{1j_r})_{\N_{1j_r}}$ is defined by the push-out
diagram $$ \xymatrix{
\Spec(k[N_{1j_1}]_{\N_{1j_r}})\ar[r]&(X_{1j_1})_{\N_{1j_r}}\\
\Spec(k[t\1N_{1j_1}]_{\N_{1j_r}})\ar[r]\ar[u]&
\Spec(k[N^-_{1j_1}]_{\N_{1j_r}})\ar[u] }
$$
while $\PP(X_{1j_r})(r)$ refers to the full subcategory of those
locally free coherent sheaves on $X_{1j_r}$ whose ranks are
non-negative multiples of $r$. The other members of the diagram
are defined similarly with respect to the corresponding polarized
monoids, and the exact functors indicated are the ones induced by
the appropriate scalar extensions, restrictions of vector bundles
to open subschemes and category embeddings.

{\it Step 5.} Fix arbitrarily a natural number $r\geq\rank N+2$.
We are going to study the lower route in the diagram
$(\ref{pdescent})_3$.

The monoids $N_{1j_r}$,$\ldots$,$N_{6j_r}$ are naturally
isomorphic to $N_1$,$\ldots$,$N_6$ so that they satisfy the
obvious analogues of the conditions (i)--(iv) in Step 2. In order
to simplify the notation we will below omit the subindex $j_r$.
Thus, the mentioned lower route looks as follows
$$
\begin{aligned}
\xymatrix{ \PP(X_1)(r)^0_{0,1}\ar[d] & &\PP(k[N_*^{\c}]_{\M})\\
\PP((X_1)_{\N_1})(r)\ar[r] &\PP((X_2)_{\N_2})(r)\ar[r]
&\PP((X_4)_{\N_4})\ar[u] }
\end{aligned}
$$
By Step 4 there is an element $\zeta_r\in
K_i(\PP(X_1)(r)^0_{0,1})$ mapping to $rx_{\M}$. Clearly,
$\zeta_r=(c_j\cdots c_{j_r})_*(z_r)$ where $z_r$ is the same as in
Step 3.

Choose a nonzero object $(P,P^-,\Theta)\in\PP(X_1)(r)^0_{0,1}$. By
Lemma \ref{matrix} and Definition \ref{diff} this object up to
isomorphism has the form
$$
P=k[N_1]^n,\ P^-=k[N_1^-]^n,\ \theta(0)\1\theta\in E(k[t\1N_1]), \
\theta(0)=
\begin{pmatrix}
t^{u_1}&0&\ldots&0\\
0&t^{u_2}&\ldots&0\\
.&.&.&.\\
0&0&\ldots&t^{u_n}
\end{pmatrix}
$$
for some $u_1,\ldots,u_n\in\{0,1\}$, where $\theta$ is the matrix
of $\Theta$ in the standard basis of $k[t\1N_1]^n$.

Since $\dim(k[t\1N_1])=\rank t\1N_1=\rank N\leq r-2\leq n-2$ the
injective stabilization estimate for $K_1$ \cite[Ch.5, \S4]{Ba}
\cite{Va} implies $\theta(0)\1\theta\in E_n(k[t\1N_1])$. It is
also clear that $\theta(0)\1\theta\in
GL_n(k[t\1N_1],\N_1k[t\1N_1]))$. Therefore, by Lemma
\ref{peculiar} there exist $\gamma^+\in GL_n(k[N_1]_{\N_1},\N_1
k[N_1]_{\N_1})$ and $\gamma^-\in GL_n(k[N^-_1]_{\N_1},\N_1
k[N^-_1]_{\N_1})$ such that $\theta(0)\1\theta=\gamma^+\gamma^-$,
the equality being considered in $GL_n(k[t\1N_1]_{\N_1})$. We have
$\theta={\big(}\theta(0)\gamma^+\theta(0)\1{\big)}
\theta(0)\gamma^-$, where by the condition (ii) in Step 2 (namely,
$t^{\pm1}(N_1(\Gamma_1)\setminus\{1\})\subset\int(N_2(\Gamma_2))$)
the inclusion $\theta(0)\gamma^+\theta(0)\1\in
GL_n(k[N_2]_{\N_2})$ holds. Therefore, by Lemma \ref{matrix} we
arrive at Claim A below.

First one notation. For $h=2,4$ we put
$$
\O_{X_h}=(k[N_h]_{\N_h},k[N_h^-h]_{\N_h},{\bf1}_{t\1k[N_h]_{\N_h}})
$$
and
$$
\O_{X_h}(-1)=(k[N_h]_{\N_h},k[N_h^-]_{\N_h},t\cdot{\bf 1}_
{t\1k[N_h]_{\N_h}}).
$$

\

{\bf Claim A.} Any object of $\PP(X_1)(r)^0_{0,1}$ is mapped to an
object of $\PP((X_2)_{\N_2})(r)$ which up to isomorphism (within
$\PP((X_2)_{\N_2})$) is of the type
$\O_{X_2}^a\oplus\O_{X_2}(-1)^b$ for some $a,b\in\ZZ_+$.

\

{\it Step 6.} Consider any morphism in $\PP((X_2)_{\N_2})$ of the
form
$$
(f^+,f^-):\O_{X_2}^a\oplus\O_{X_2}(-1)^b\to
\O_{X_2}^c\oplus\O_{X_2}(-1)^d.
$$
We let $(f^+_{ij})$ and $(f^-_{ij})$ denote the matrices of $f^+$
and $f^-$ in the corresponding standard bases. Then we see that
there are systems $u_1,\ldots,u_{a+b}\in\{0,1\}$ and
$v_1,\ldots,u_{c+d} \in\{0,1\}$ such that $$
f_{ij}^+t^{v_j}=f_{ij}^-t^{u_i},\ i\in[1,a+b],\ j\in[1,c+d]. $$
(Compare with the equations (\ref{poldesc})$_1$ and
(\ref{poldesc})$_2$.) By the condition (ii) in Step 2 (namely,
$t^{\pm1}(N_2(\Gamma_2)\setminus\{1\})\subset\int(N_3(\Gamma_3))$)
we get the inclusions
$$
f_{ij}^+\in k[N_3(\Gamma_3)]_{\N_3}+k\cdot t\ \ \text{and}\ \
f_{ij}^-\in k[N_3^-(\Gamma_3)]_{\N_3}+k\cdot
t\1\eqno(\ref{pdescent})_4
$$
for all $i$ and $j$.

\

{\bf Claim B.} Any short exact sequence in $\PP((X_2)_{\N_2})$ of
the type
$$
0\to\O_{X_2}^a\oplus\O_{X_2}(-1)^b
\to\O_{X_2}^{a'}\oplus\O_{X_2}(-1)^{b'}
\xto{f}\O_{X_2}^{a''}\oplus\O_{X_2}(-1)^{b''}\to0
$$
is mapped to a split short exact sequence in $\PP((X_4)_{\N_4})$.

\

By Lemma \ref{split} any such a sequence reduces modulo $\N_2$ to
a split sequence in $\PP(\Pp^1_k)_{0,1}$. Let
$\pi:\PP((X_2)_{\N_2})\to\PP(\Pp^1_k)$ denote the exact functor,
obtained by reduction modulo $\N_2$, and let $g$ be a monomorphism
that splits $\pi(f)$.

We have the scalar extension functor
$\iota:\PP(\Pp^1_k)\to\PP((X_2)_{\N_2})$. The composite
$\pi\circ\iota$ is isomorphic to the identity functor on
$\PP(\Pp^1_k)$. We can assume
$$
\pi(f\circ\iota(g))={\bf
1}_{\O^{a_3}\oplus\O(-1)^{b_3}}.\eqno(\ref{pdescent})_5
$$
Claim B is proved once we show that the composite $f\circ\iota(g)$
is mapped to an automorphism of
$\O_{X_4}^{a''}\oplus\O_{X_4}(-1)^{b''}$ under the functor
$\PP((X_2)_{\N_2})\to\PP((X_4)_{\N_4})$.

It is easily observed that the positive and negative components of
$\iota(g)$ have correspondingly entries from $k+k\cdot t$ and
$k+k\cdot t^{-1}$. But it then follows from the condition (ii) in
Step 2 (namely,
$t^{\pm1}(N_3(\Gamma_3)\setminus\{1\})\subset\int(N_4(\Gamma_4))$),
the inclusion $(\ref{pdescent})_4$ and the equality
$(\ref{pdescent})_5$ that both positive and negative components of
$f\circ\iota(g)$ are defined over the local ring
$k[N_4(\Gamma_4)]_{\N_4}=k[N_4^-(\Gamma_4)]_{\N_4}$, that is the
corresponding matrices w.r.t. the standard bases have entries from
this local ring. Since these matrices are invertible modulo $\N_4$
they are invertible themselves.

{\it Step 7.} Let $\E\subset\PP((X_4)_{\N_4})$ denote the full
subcategory whose objects are $\O_{X_4}^a\oplus\O(-1)_{X_4}^b$,
$a,b\in\ZZ_+$ (in the selfexplanatory notation). This is an
additive category. We equip it with an exact structure with
respect to the class of split short sequences -- any additive
category carries such an exact structure, as shown by the Yoneda
embedding ($*\mapsto\Hom(-,*)$) in the category of contravariant
additive functors with values in abelian groups. By the  previous
steps we have the natural commutative triangle of exact categories
and exact functors
$$
\begin{aligned}
\xymatrix{
&\E\ar[d]\\
\PP(X_1)(r)^0_{0,1}\ar[ur]\ar[r]&\PP((X_4)_{\N_4}) }
\end{aligned}
\eqno(\ref{pdescent})_6
$$
Consider a morphism
$f=(f^+,f^-):\O_{X_4}^{a_1}\oplus\O(-1)_{X_4}^{b_1}\to
\O_{X_4}^{a_2}\oplus\O(-1)_{X_4}^{b_2}$. For the corresponding
matrices we have the equality
$$
\begin{pmatrix}
f^+_{a_1\times a_2}&tf^+_{a_1\times b_2}\\
f^+_{b_1\times a_2}&tf^+_{b_1\times b_2}
\end{pmatrix}=
\begin{pmatrix}
f^-_{a_1\times a_2}&f^-_{a_1\times b_2}\\
tf^-_{b_1\times a_2}&tf^-_{b_1\times b_2}
\end{pmatrix}.
$$
Since $N_4$ is polarized (and
$U(k[N_4(\Gamma_4)]_{\N_4})=k^*+\N_4k[N_4(\Gamma_4)]_{\N_4}$) we
get:
\begin{itemize}
\item
the entries of $f^+_{a_1\times a_2}$ and $f^+_{b_1\times b_2}$
belong to $k[N_4(\Gamma_4)]_{\N_4}$,
\item
the entries of $f^+_{a_1\times b_2}$ belong to
$\N_4k[N_4(\Gamma_4)]_{\N_4}\cap t\1\N_4k[N_4(\Gamma_4)]_{\N_4}$,
\item
the entries of $f^+_{b_1\times a_2}$ belong to
$k+\N_4k[N_4(\Gamma_4)]_{\N_4}+t\N_4k[N_4(\Gamma_4)]_{\N_4}+k\cdot
t$.
\end{itemize}
Conversely, any matrix over $k[N_4]_{\N_4}$ of size
$(a_1+b_1)\times(a_2+b_2)$, whose entries satisfy the three
conditions above, is a positive component of a unique morphism in
$\E$.

Consider the subring $\Lambda_0=\{(\phi_{uv})\}\subset
M_{2\times2}(k[N_4]_{\N_4})$ where
\begin{itemize}
\item
$\phi_{11},\phi_{22}\in k[N_4(\Gamma_4)]_{\N_4}$,
\item
$\phi_{12}\in\N_4k[N_4(\Gamma_4)]_{\N_4}\cap
t\1\N_4k[N_4(\Gamma_4)]_{\N_4}$,
\item
$\phi_{21}\in
k+\N_4k[N_4(\Gamma_4)]_{\N_4}+t\N_4k[N_4(\Gamma_4)]_{\N_4}+k\cdot
t$.
\end{itemize}

Consider the full subcategory $\E_{\Lambda_0}\subset\E$ of objects
of the type $O_{X_4}^a\oplus\O(-1)_{X_4}^a$. By the mentioned
description of morphisms in $\E$ we have the natural equivalence
of categories $\E_{\Lambda_0}\approx{\Bbb F}(\Lambda_0)$ -- the
category of left free $\Lambda_0$-modules (after thinking of
elements of $\Lambda_0$ as the corresponding endomorphisms of
$O_{X_4}\oplus\O(-1)_{X_4}$).

Clearly, $\E_{\Lambda_0}$ is a cofinal subcategory of the exact
category $\E$ in which all exact sequences split by definition.
Therefore, by \cite[\S The plus construction]{Gra} we have
$K_i(\E_{\Lambda_0})=K_i(\E)$. (Recall $i\geq1$; for the
Grothendieck groups we have $K_0(\E)=\ZZ\oplus\ZZ$ and
$K_0(\E_{\Lambda_0})=\ZZ$.)

We have the exact functor
$\PP(\E_{\Lambda_0})\to\PP(k[N^{\c}_*]_{\M})$ given by picking up
the ``positive'' part in the objects and morphisms in
$\E_{\Lambda_0}$. In other words, we mean the composite functor
$\E_{\Lambda_0}\subset\E\to\PP((X_4)_{\N_4})\to\PP(k[N_4]_{\N_4})
\to\PP(R[N^{\c}_*]_{\M})$, where the third functor is given by
restriction to an open subscheme and the fourth functor
corresponds to a scalar extension.

{\it Step 8.} Consider the two subrings
$\Lambda_4=\Lambda_{(t,\Gamma_4,N_4)}$ and $\Lambda_5=
\Lambda_{(t,\Gamma_5,N_5)}$ of the matrix ring
$M_{2\times2}(k[N_*^{\c}])$. Put
$S=k^*+k{\big(}\int(N_5(\Gamma_5){\big)}\subset
k[N_5(\int(\Gamma_5))]$. We regard $S$ a central multiplicative
set in $M_{2\times2}(k[N_*^{\c}])$ through the diagonal embedding.

It is easily observed that $(\Lambda_4)_{\N_4}=\Lambda_0$. (For an
essentially equivalent equality see \cite[Proposition 2.14, Step
5]{Gu2}.) The inclusion $\Gamma_4\subset\int(\Gamma_5)$ implies
$\Lambda_0\subset S\1\Lambda_5$. Therefore, we have the
commutative diagram (with obvious morphisms)
$$
\xymatrix{
&\Lambda_5\ar[rrr]\ar[ddd]&&&M_{2\times2}(k[N_5])\ar[ddd]\\
&&\Lambda_4\ar[ul]\ar[r]\ar[d]&M_{2\times2}(k[N_4])\ar[ur]\ar[d]&\\
&&\Lambda_0\ar[dl]\ar[r]&M_{2\times2}(k[N_4])_{\N_4}\ar[dr]&\\
&S\1\Lambda_5\ar[rrr]&&&S\1M_{2\times2}(k[N_5])
}\eqno(\ref{pdescent})_7
$$
By Lemmas \ref{int} and \ref{kar} the boundary of
$(\ref{pdescent})_7$, which consists of
$k[N_5(\int(\Gamma_5))]$-algebra homomorphisms, is a Karoubi
square. Hence the long exact sequence
$$
\cdots\to K_i(\Lambda_5)\to K_i(k[N_5])\oplus K_i(S\1\Lambda_5)\to
K_i(S\1k[N_5])\to\cdots\eqno(\ref{pdescent})_8
$$
By Step 7 we have the natural commutative diagram
$$
\xymatrix{
&K_i(k[N_4])\ar[rd]&&\\
K_i{\big(}\PP(X_1)(r)_{0,1}^0{\big)}\ar[ru]\ar[rd]&&
K_i(k[N_4]_{\N_4})\ar[r]&K_i(k[N_*^{\c}]_{\M})\\
&K_i(\Lambda_0)\ar[ru]&& }
$$
Also, using $\Gamma_5\subset\int(\Phi(M))$, we have the
homomorphisms with the same composite
$$
K_i(\PP(X_1)(r)^0_{0,1})\to K_i(k[N_5])\to K_i(S\1k[N_5])\to
K_i(k[N_*^{\c}]_{\M}).
$$
Therefore, using the element $\zeta_r\in
K_i{\big(}\PP(X_1)(r)_{0,1}^0{\big)}$ (see Step 5) and the
sequence $(\ref{pdescent})_8$, we see that
$rx\in\I{\big(}K_i(\Lambda_5)\to K_i(k[N_*^{\c}]){\big)}$ for
$r\geq\rank N+2$, i.~e. $x\in\I{\big(}K_i(\Lambda_5)\to
K_i(k[N_*^{\c}]){\big)}$.

{\it Step 9.} For a natural number $c$ put
$\Lambda^{c}=\Lambda_{(t^c,\Gamma_5,N_5^c)}$. Clearly,
$N_5^{c_1\cdots c_j}\subset N_*$ for $j\gg0$. Then Step 8 implies
$(c_1\cdots c_j)_*(z)\in\I{\big(}K_i(\Lambda_5^{c_1\cdots c_j})\to
K_i(k[N_*]){\big)}$ for $j\gg0$. It is easily observed that all
the polarized monoids that have shown up above work equally well
for the elements $c_*(z)\in K_i(k[N_*])$ for {\it any} natural
number $c$. Therefore, for all $j\gg0$ and all $j'\geq0$ we have
$(c_1\cdots c_{j+j'})_*(z)\in\I{\big(}K_i(\Lambda_5^{c_1\cdots
c_j})\to K_i(k[N_*]){\big)}$. By the conditions (i) and (ii) in
Step 2 (namely, $\Gamma_6\subset\int(\Phi(M))$ and
$t^{\pm1}(N_5(\Gamma_5)\setminus\{1\})\subset\int(N_6(\Gamma_6))$)
we have the inclusion
$$
\Lambda'_{(t^{c_1\cdots c_j},\Gamma_5,N_5^{c_1\cdots c_j})}\subset
M_{2\times2}(k[M_*]).
$$
Therefore, the polarized monoid $(t^{c_1\cdots
c_j},\Gamma_5,N_5^{c_1\cdots c_j})$ is the desired one for
$j\gg0$.
\end{proof}

\section{The action of Witt vectors}\label{WITT}

Assume $R$ is a general commutative ring. Stienstra \cite{St} has
studied a continuous $\W(R)$-module structure on $NK_i(R)$,
$i\in\ZZ_+$. (Such actions more or less implicitly were previously
defined by Bloch \cite{Bl}.) Recall the additive group of $\W(R)$
is the multiplicative group $1+TR[[T]]$ and the multiplicative
structure is determined by
$$
(1-rT^m)\star(1-sT^n)=(1-r^{n/d}s^{m/d}T^{mn/d})^d, \quad r,s\in
R,\quad d=\gcd(m,n).
$$
We have the decreasing ideal filtration $I_m(R)=1+T^mR[[T]]$.
``Continuous'' here means that the annihilator of any element
$z\in NK_i(R)$ contains $I_m(R)$ for some $m$.

Weibel \cite{W2} has generalized these operations to the graded
situation as follows. Assume $A=A_0\oplus A_1\oplus\cdots$ is a
graded, not necessarily commutative ring and $R\subset A_0$ is a
subring in the center of $A$. Then there is a {\it functorial}
continuous $\W(R)$-module structure on $K_i(A,A^+)$. (Here
$A^+=0\oplus A_1\oplus A_2\oplus\cdots$.)

In the special case $A=A_0[T]$ the action of $1-rT^n\in\W(R)$,
$r\in R$ on $NK_i(A_0)$ is the effect of the composite functor
$$
\Pp(A_0[T])\xto{t_n}\Pp(A_0[T])\xto{r_*}\Pp(A_0[T])\xto{(v_n)_*}\Pp(A_0[T])
$$
where $v_n:A_0[T]\to A_0[T]$ is given by $T\mapsto T^n$,
$r:A_0[T]\to A_0[T]$ is given by $T\mapsto rT$ and $t_n$ is the
scalar restriction through $v_n$. This determines the action of
the whole $\W(R)$ because any element $\omega(T)\in\W(R)$ has a
unique convergent expansion $\omega(T)=\Pi_{n\geq1}(1-r_nT^n)$.
Moreover, since for $\omega(T)\in I_m(R)$ the expansion has the
form $\omega(T)= \Pi_{n\geq m}(1-r_nT^n)$, we conclude that
$$
I_m(R)NK_i(A_0)\subset\I{\big(}K_i(A_0[\{T^n\}_{n\geq m}]\to
K_i(A_0[T]){\big)},\quad m\in\NN.\eqno{(\ref{WITT})_1}
$$
Now assume a natural number $n$ is invertible in $A_0$ and $z\in
NK_i(A_0)$ is in the image of $NK_i((v_n)_*)$. Since $n$ is
invertible also in $\W(R)$ we have $\frac1nz\in NK_i(A_0)$. On the
other hand the composite $NK_i(t_n)\circ
NK_i((v_n)_*):NK_i(A_0[T]\to NK_i(A_0[T])$ is the multiplication
by $n$. Therefore, $z\in\I{\big(}NK_i((v_n)_*)\circ NK_i(t_n)\circ
NK_i((v_n)_*){\big)}\subset\I{\big(}NK_i((v_n)_*)\circ
NK_i(t_n){\big)}=(1-T^n)\star NK_i(A_0)$, $\star$ referring to the
$\W(R)$-action. In particular, if $\QQ\subset A_0$ we have the
equality
$$
I_m(R)NK_i(A_0)=\I{\big(}K_i(A_0[\{T^n\}_{n\geq m}])/K_i(A_0)\to
NK_i(A_0){\big)},\quad m\in\NN.\eqno{(\ref{WITT})_2}
$$
Assume for simplicity $R=k$ -- a characteristic 0 field. Since
$\W(k)\approx\Pi_1^{\infty}k$ ({\it ghost map}) with the product
topology, a continuous $\W(k)$-module is just a $k$-vector space
with a grading by natural numbers. Further, the {\it Weibel map}
of graded algebras $\w:A\to A[T]$, $\sum_ja_j=\sum_ja_jT^j$
($A[T]$ being graded w.r.t. the degrees of $T$) induces an {\it
embedding} of $\W(k)$-modules $K_i(A,A^+)\to NK_i(A)$ -- a
consequence of the fact that $\w$ splits the non-graded ring
homomorphism $A[T]\to A$, $T\mapsto1$. Using $(\ref{WITT})_2$ (for
the ring $A$) we obtain the embedding $K_i(A_{[n]},A_{[n]}^+)\to
I_n(k)NK_i(A)$ where $A_{[n]}=A_0\oplus0\oplus\cdots\oplus0\oplus
A_{n}\oplus A_{n+1}\oplus\cdots$. Summarizing, we get

\begin{proposition}\label{totaro}
The following hold in the category of $k$-vector spaces:
\begin{itemize}
\item[(a)]
$K_i(A,A^+)=\bigoplus_{j=1}^{\infty} V_j$,
\item[(b)]
$I_n(k)K_i(A,A^+)=0\oplus\cdots\oplus0\oplus V_n\oplus
V_{n+1}\oplus\cdots$, \ \ $n\in\NN$,
\item[(c)]
$\I{\big(}K_i(A_{[n]},A_{[n]}^+)\to K_i(A,A^+{\big)}\subset
I_n(k)K_i(A,A^+), \quad n\in\NN$.
\end{itemize}
\end{proposition}

We would like to have a weaker version of $(\ref{WITT})_1$ in the
graded situation as follows.

\begin{question}\label{witt}
Let $k$ be a characteristic 0 field and $A=A_0\oplus
A_1\oplus\cdots$ be a graded not necessarily commutative
$k$-algebra. Assume $A_1=A_2=\cdots= A_{m-1}=0$, $\dim_k(A_m)=1$
for some $m\in\NN$, and $\dim_k(A_j)<\infty$ for all $j\geq0$.
Does one have the inclusions
$I_n(k)K_i(A,A^+)\subset\I{\big(}K_i(A_{[m+1]})\to K_i(A){\big)}$
for $n\gg0$?
\end{question}

The relevance of Question \ref{witt} is that the positive answer
to it would complete the proof of Conjecture \ref{conj} for fields
of characteristic 0.

To see this consider a pyramidal extension of monoids $M\subset N$
and an element $z\in K_i(k[N_*])$ with the property $z(0)\in
K_i(k)$. Let $(t,\Gamma,L)$ be as in Theorem \ref{pdescent}. We
know that there is a grading $k[N_*]=k\oplus R_1\oplus\cdots$ such
that $N_*$ consists of homogeneous elements. It is easily observed
that we can choose the grading of $k[N_*]$ in such a way that the
pole $t$ is the element of the smallest positive degree within the
polarized monoid $L$. By Proposition \ref{totaro} we have
$c_*(z)\in I_c(k)K_i(k[N_*],k[N_*]^+)$ for $c\in\NN$. Further, we
have the induced gradings
$\Lambda_{(t,\Gamma,L)}=\Lambda_0\oplus\Lambda_1\oplus\cdots$ and
$\Lambda'_{(t,\Gamma,L)}=\Lambda'_0\oplus\Lambda'_1\oplus\cdots$
and these rings are in the same relation as $A$ and $A_{[m+1]}$ in
Question \ref{witt}. Therefore, the positive answer to this
question (and the Morita invariance of the $\W(k)$ actions,
\cite{W2}) would imply $(c_1\cdots
c_j)_*(z)\in\I{\big(}K_i(\Lambda'_{(t,\Gamma,L)})\to
K_i(k[N_*]){\big)}$ for $j\gg0$. But the mentioned map factors
through $K_i(k[M_*])$ and, hence, the pyramidal descent is
achieved.

In one situation, which is to some extent dual to the situation
considered so far, we have the following partial result on
Question \ref{witt}. It will be used in \S\ref{BP} -- in a
different way but again with essential use of Theorem
\ref{pdescent} -- to prove Conjecture \ref{conj} for certain class
of non-simplicial monoids.

\begin{lemma}\label{findeg}
Let $k$ be a number field (i.~e. a finite extension of $\QQ$) and
$A=A_0\oplus\ A_1\oplus\cdots$ be a graded $k$-algebra such that
$\dim_kA<\infty$. Then $K_i(A,A^+)$ is a finite dimensional
$k$-vector space. In particular, $I_n(k)K_i(A,A^+)=0$ for $n\gg0$.
\end{lemma}

\begin{proof} $A^+$ is a nilpotent ideal. So by Goodwillie's
result \cite{Go} we have the isomorphism of abelian groups
$K_i(A,A^+)\approx HC_{i-1}(A,A^+)$, the right hand side being a
finite dimensional $\QQ$-vector space.
\end{proof}

\begin{remark}\label{remark}\label{kabi}
(a) It is very natural to expect that the answer to Question
\ref{witt} is `yes' when $k$ is a number field. This would proof
Conjecture \ref{conj} for toric varieties over such fields.
Notice, however, that one needs to modify the argument in Step 1
in the proof of Theorem \ref{pdescent} so that the inductive step
no longer increases the transcendence degree of the ground field.
Here we only mention that this is always possible. For the special
class of bipyramidal cones the details are included in Step 1 in
the proof of Theorem \ref{bp} below.

(b) We also expect that  $I_n(k)K_i(A,A^+)=0$ for $n\gg0$ provided
$\chara k=0$ and $\dim_kA<\infty$. This would extend Theorem
\ref{bp} to all characteristic $0$ fields.
\end{remark}

\section{Bipyramidal monoids}\label{BP}

Let $P\subset\RR^r$ be a finite convex polytope of dimension $\dim
P<r$ and $\sigma\subset\RR^r$ be a closed segment (a 1-dimensional
polytope). Assume $\int(\sigma)\cap\int(P)$ is a point. Then the
polytope $Q=\conv(\sigma\cup P)$ will be called a {\it bipyramid}
over $P$. The class of polytopes $\BP$ is defined recursively by
the condition:
\begin{itemize}
\item
Points are in $\BP$ and if $P\in\BP$ then pyramids and bipyramids
over $P$ are in $\BP$.
\end{itemize}
Every positive dimensional polytope $P\in\BP$ admits a finite, not
necessarily unique sequence of polytopes $P_1\subset
P_2\subset\cdots\subset P_d=P$, $d=\dim P$ such that $P_1$ is a
segment and $P_s$ is either a pyramid or a bipyramid over
$P_{s-1}$,  $s\in[2,\dim P]$. To such a sequence we associate a
sequence $\sigma(P)$ of length $d-1$, consisting of 0s and 1s, as
follows: the $s$th member of $\sigma(P)$ is 1 if $P_{d-s+1}$ is a
bipyramid over $P_{d-s}$ and 0 otherwise. For instance, simplices
are characterized by the condition $\sigma(P)=0\ldots0$,
$\sigma(P)=11$ means that $P$ is an octahedron, $\sigma(P)=01$ if
$P$ is a pyramid over a square etc. The sequence $\sigma(P)$ is
{\it uniquely} determined by the combinatorial type of $P$. To see
this notice that no bipyramid can simultaneously be a pyramid
because arbitrary facet of a bipyramid admits at leat two vertices
not in this facet. In particular, the first member of $\sigma(P)$
is uniquely determined. But then we can apply induction on the
length of such sequences because the combinatorial type of the
preceding member $P_{d-1}$ is determined by that of $P$ -- an easy
observation. As a result, we can identify $\sigma(P)$ with the
combinatorial type of $P$. Therefore, in dimension $r>0$ we have
$2^{r-1}$ different $\BP$-combinatorial types  of which one is
simplicial.

We order the set of  $\BP$-combinatorial types with respect first
to the dimension and then to the lexicographic order of the
$\sigma$-sequences.

A finite polyhedral pointed cone $C\subset\RR^r$ will be called
{\it of type $\BP$} if its polytopal cross section is in $\BP$. A
monoid $N$ without nontrivial units will be called {\it of type
$\BP$} if the polytope $\Phi(N)$ is in $\BP$. Obviously, each
$\BP$-combinatorial type admits infinitely many nonisomorphic
monoids, whose $\Phi$-polytopes are of this combinatorial type.

\begin{theorem}\label{bp}
Let $\QQ\subset k$ be an algebraic extension of fields, $N$ be a
monoid of type $\BP$ and $i\in\NN$. Then $\NN$ acts nilpotently on
$K_i(k[N])$. If, in addition, $N$ is normal and finitely generated
then $\NN$ acts nilpotently on $K_i(\Proj(k[N]))$ for arbitrary
grading $k[N]=k\oplus R_1\oplus\cdots$ making the elements of $N$
homogeneous.
\end{theorem}
\noindent Here the nilpotence of the action in the projective case
is understood in the sense of Proposition \ref{global}. Also, we
remark that the normality assumption on $N$ in the second half of
Theorem \ref{bp} can actually be dropped.

Among the toric varieties over $k$ covered by this theorem, the
cone over the Segre embedding $\Pp_k^1\times\Pp^1_k\to\Pp^3_k$ is
the simplest which is not simplicial and
$$
\Proj(k[U,V,W,X,Y,Z]/(UX-WZ,UX-VY))\quad
(\deg(U)=\cdots=\deg(Z)=1)
$$
is the simplest projective variety whose standard affine covering
involves no simplicial toric variety.

We need a preparation. Let $C\subset\RR^r$ be a rational
bipyramidal cone and $C=C_1\cup C_2$ be a decomposition into
pyramidal cones. (We call a cone `bipyramidal' or `pyramidal' if
its polytopal cross section is such.) Consider a nonzero element
$t\in C_2\cap\ZZ^r$ such that $\RR_+t$ is the unique extremal ray
of $C_2$ not in $C_1\cap C_2$. Denote by $\mathcal L$ the monoid
$C_1\cap\ZZ^r$ and by $l$ the extremal ray of $C_1$ not in
$C_1\cap C_2$.

Consider the subring $\Lambda_{t,C}=\{(\phi_{uv})\}\subset
M_{2\times2}(k[\ZZ^r])$  given by
\begin{itemize}
\item
$\phi_{11},\phi_{22}\in k[{\mathcal L}]$,
\item
$\phi_{12}\in k[{\mathcal L}]\cap t\1k[{\mathcal L}]$,
\item
$\phi_{21}\in k[{\mathcal L}]+t k[{\mathcal L}]$.
\end{itemize}
Assume the following holds
\begin{itemize}
\item
$\omega=l\cap{\big(}-t+(C_1\cap C_2){\big)}\in\ZZ^r$.
\end{itemize}
In this situation $\omega$ and $t+\omega$ both belong to
${\mathcal L}$. We also have  $k[{\mathcal L}]\cap t\1k[{\mathcal
L}]=\omega k[{\mathcal L}]$ and $c\cdot
t+(c-1)\cdot\omega=t+(c-1)\cdot(t+\omega)\in t+{\mathcal L}$ for
$c\in\NN$. It follows that for every natural number $c$ we have
the $k$-algebra endomorphism:
$$
\tilde c:\Lambda_{t,C}\to\Lambda_{t,C},\quad (\phi_{uv})\mapsto
\begin{pmatrix}
c_{\#}(\phi_{11})&\omega^{-c+1}c_{\#}(\phi_{12})\\
&\\
\omega^{c-1}c_{\#}(\phi_{21})&c_{\#}(\phi_{22})
\end{pmatrix}
$$
where $c_{\#}=k[-^c]$. In this notation we have the following

\begin{lemma}\label{endo}
Let $C'\subset\RR^r$ be a polyhedral rational pointed cone such
that $C_1\setminus\{0\}\subset\int(C')$. Then $\I(\tilde c)\subset
M_{2\times2}(k[(C'\cap\ZZ^r])_*]$ for $c\gg0$.
\end{lemma}

\begin{proof}
Since $t+\omega\in\int(C')$ we have $\omega^{c-1}t^c\in\int(C')$
for $c\gg0$ and this yields the desired inclusion.
\end{proof}

Let $N$ be arbitrary finitely generated, normal monoid with
trivial $U(N)$ and $C(N)=C'\cup C''$ be a decomposition into two
nondegenerate pyramidal cones, sharing a base facet. Put
$M=N(C')$. Then $M\subset N$ is a pyramidal extension of monoids.
Assume $(t,\Gamma,L)$ is a polarized monoid such that
$L\setminus\{0\}\subset\int(N)$, $\Gamma\subset\int(\Phi(M))$, and
$t\in C''\setminus C'$. As usual, $\c=(c_1,c_2,\ldots)$ where
$c_1,c_2,\ldots\geq2$.

\begin{lemma}\label{bpappro}
There exists a rational bipyramidal cone $C$ with the property
$C\setminus\{0\}\subset\int(C(N))$ and admitting a decomposition
$C=C_1\cup C_2$ into pyramidal cones in such a way that the
following hold (in the real space $\RR\otimes\gp(N))$:
\begin{itemize}
\item[(a)]
$\Gamma\subset\int(C_1)$, $C_1\setminus\{0\}\subset\int(C')$, and
$L\subset C$,
\item[(b)]
$C_2=\RR_+t+(C_1\cap C_2)$,
\item[(c)]
$l\cap{\big(}-t+(C_1\cap C_2){\big)}\in\gp(N)^{\c}$ where
$l\subset C_1$ is the extremal ray not in $C_1\cap C_2$.
\end{itemize}
\end{lemma}

\begin{proof}
We can gradually approximate $\int(C')$ by appropriate pyramidal
rational cones $C_\alpha$, $\alpha\in\NN$ such that
$C_\alpha\setminus\{0\}\subset\int(C')$. Since
$\gp(N)^{\c}\subset\RR\otimes\gp(N)$ is a {\it dense} subset, we
can also keep the condition
$l_\alpha\cap(-t+F_\alpha)\in\gp(N)^{\c}$ satisfied, where
$F_\alpha\subset C_\alpha$ is the facet `close' to $C'\cap C''$
and $l_\alpha$ is the extremal ray of $C_\alpha$ not in
$F_\alpha$. Then $C=C_\alpha+\RR_+t$ is the desired bipyramidal
cone for $\alpha\gg0$ whose desired decomposition into pyramidal
cones is $C=C_\alpha\cup(F_\alpha+\RR_+t)$
\end{proof}

We also need several results on triangular matrix rings.

Berrick and Keating showed \cite{BeKe} that for not necessarily
commutative rings $A$ and $B$ and an $A - B$ bimodule $U$ the
embedding of rings
$$
A\oplus B\to
\begin{pmatrix}
A&0\\
U&B
\end{pmatrix}
$$
induces isomorphisms on $K$-groups. Further, let $R$ be a
commutative ring and $I\subset R$ be a {\it projective} ideal. Put
$$
A=\begin{pmatrix}
R&I\\
R&R
\end{pmatrix}\quad\text{and}\quad B=
\begin{pmatrix}
R/I&0\\
R/I&R/I
\end{pmatrix}.
$$
For every natural index $i$ we have the group homomorphisms
\begin{itemize}
\item
$f:K_i(A)\to K_i(M_{2\times2}(R))=K_i(R)$,
\item
$g:K_i(A)\to K_i(B)$,
\item
$h:K_i(B)\to K_i(M_{2\times2}(R/I))=K_i(R/I)$.
\end{itemize}

Next is the commutative case of Keating's result.

\begin{theorem}[Keating \cite{Ke}]\label{keat}
For any natural index $i$ we have
\begin{itemize}
\item[(a)]
$h{\big(}(a,b){\big)}=a+b$ for $a,b\in K_i(R/I)$. In particular,
$\Ker(h)\approx K_i(R/I)$.
\item[(b)]
$f:K_i(A)\to K_i(R)$ is a split epimorphism.
\item[(c)]
$g$ restricts to an isomorphism $\Ker(f)\to\Ker(h)$. In
particular, $K_i(A)=K_i(R)\oplus K_i(R/I)$.
\end{itemize}
\end{theorem}

{\it Proof of Theorem \ref{bp}}. There is no loss of generality in
assuming that $\QQ\subset k$ is a finite extension.

{\it Step 1.} First we consider the affine case.

We have well ordered the set of combinatorial types of the cones
at the issue and the proof can be carried out by induction with
respect to this order. Theorem \ref{simplok} gives the result for
monoids $N$ such that $\sigma(\Phi(N))=0\ldots0$. This includes
the case $\rank N\leq2$.

All faces of a polytope in $\BP$ are in $\BP$. Therefore, by the
obvious version of Lemma \ref{inter} for monoids of type $\BP$ it
suffices to show that $\NN$ acts nilpotently on $K_i(k[N_*])$
where $N$ is in addition a finitely generated and normal monoid.

In order to be able to apply Theorem \ref{pdescent} we have to
resolve one difficulty -- the induction hypothesis on $\rank N$,
used in Step 2 in the proof of Theorem \ref{pdescent}, {\it
increases} the transcendence degree of the ground field. Here we
develop an alternative, circumventing such an increase for our
special cones.

There is no loss of generality in assuming
$\Phi(N)=0\ldots01\delta_s\ldots\delta_{d-1}$, for some
$s\in[1,d]$ ($d=\rank(N)-1$). The cone $C(N)$ decomposes into two
pyramidal cones $C(N)=C'\cup C''$ whose polytopal cross sections
are both of type $0\ldots00\delta_s\ldots\delta_{d-1}$. Let
$M\subset N$ be the corresponding pyramidal extension so that
$C(M)=C'$. By the obvious adaptation of Theorem \ref{approxB} to
the monoids of type $\BP$ we can assume that the polarized monoids
$(t_\alpha,\Gamma_\alpha,N_\alpha)$ are such that the
$\Gamma_\alpha$ are pyramids whose bases approximate $C'\cap C''$
and $\sigma(\Gamma_\alpha)=0\ldots00\delta_s\ldots\delta_{d-1}$,
notation as in Theorem \ref{approxB}. In this situation the
monoids $(t_\alpha^{-1}N_\alpha)\setminus(t_\alpha^{-1})^\NN$ are
filtered unions of finitely generated normal monoids
$N_{\alpha\beta}$ of type $\BP$ such that
$\sigma(\Phi(N_{\alpha\beta}))=\sigma(\Gamma_\alpha)=
0\ldots00\delta_s\ldots\delta_{d-1}<\sigma(\Phi(N))$. Therefore,
for each element $z\in K_i(k[N_*])$, $z(0)=0$, the condition (iv)
in Step 2 in the proof of Theorem \ref{pdescent} is achieved by
the induction assumption.

All the other arguments in the proof of Theorem \ref{pdescent} go
through with respect to the system
$(t_\alpha,\Gamma_\alpha,N_\alpha)$. So we reach the situation
when:
\begin{itemize}
\item
$N$ is a finitely generated normal monoid without nontrivial
units,
\item
$C(N)=C'\cup C''$ is a decomposition into two pyramidal cones
sharing a base facet and $\sigma(\Phi(M))<\sigma(\Phi(N))$ where
$M=C'\cap N$,
\item
$(t,\Gamma,L)$ is polarized monoid such that $L\subset N_*$ and
$\Gamma\subset\int(\Phi(M))$,
\item
$z\in K_i(k[N_*])$, $z(0)=0$, and $(c_1\cdots
c_j)_*(z)\in\I{\big(}K_i(\Lambda_{(t,\Gamma,L)}\to
K_i(k[N_*]){\big)}$ for $j\gg0$ (notation as in Theorem
\ref{pdescent}).
\end{itemize}
By induction assumption it suffices to achieve the pyramidal
descent for the  extension $M\subset N$, i.~e. $(c_1\cdots
c_j)_*(z)\in\I{\big(}K_i(k[M_*]\to k[N_*]{\big)}$ for $j\gg0$.

{\it Step 2.} Fix $j_1$ such that $(c_1\cdots
c_{j_1})_*(z)\in\I{\big(}K_i(\Lambda_{(t,\Gamma,L)}\to
K_i(k[N_*]){\big)}$. Next we choose a bipyramidal cone
$C\setminus\{0\}\subset\int(C(N))$ as in Lemma \ref{bpappro} with
respect to the cone $C'$, the polarized monoid $(t,\Gamma,L)$, and
the sequence $\c'=(c_{j_1+1},c_{j_1+2},\ldots)$. Put
$\omega=l\cap{\big(}-t+(C_1\cap C_2){\big)}$ (notation as in that
lemma) and choose $j_2>j_1$ such that $\omega\in(c_{j_1+1}\cdots
c_{j_2})\1\gp(N)$. For simplicity of notation put
$\kappa=c_{j_1+1}\cdots c_{j_2}$ and consider the ring
$\Lambda_{t^\kappa,C}$ as in Lemma \ref{endo} with respect to the
lattice $\gp(N)\approx\ZZ^{\rank N}$. We have the diagram of ring
embeddings
$$
\xymatrix{ \Lambda_{(t^\kappa,\Gamma,L^\kappa)}
\ar[r]&\Lambda_{t^\kappa,C}\ar[rd]&&&\\
&&A\ar[r]&M_{2\times2}(k[{\mathcal N}])\ar[r]&M_{2\times2}(k[N_*])\ ,\\
&k[{\mathcal N}]\ar[ru]_{\Delta}&&& }
$$
where $\mathcal N=C\cap\gp(N)$,
$$
A=
\begin{pmatrix}
k[{\mathcal N}]&\omega^\kappa k[{\mathcal N}]\\
\\
k[{\mathcal N}]&k[{\mathcal N}]
\end{pmatrix}\subset M_{2\times2}(k[{\mathcal N}]),\quad \Delta(a)=
\begin{pmatrix}
a&0\\
0&a
\end{pmatrix},
$$
and all the other homomorphisms are the indentity embeddings. As
observed in Step 9 in the proof of Theorem \ref{pdescent},
$(c_1\cdots
c_{j_1}c)_*(z)\in\I{\big(}K_i(\Lambda_{(t,\Gamma,L)}\to
(k[N_*]){\big)}$ for any natural number $c$. By choosing
$c=c_{j_2+1}c_{j_2+2}\cdots c_j$ for $j>j_2$ we get
$$
(c_1\cdots c_{j_1}c_{j_1+1}\cdots c_{j_2}c_{j_2+2}\cdots
c_j)_*(z)\in\I{\big(}K_i(\Lambda_{(t^\kappa,\Gamma,L^\kappa)} \to
M_{2\times2}(k[N_*]){\big)}
$$
for all $j>j_2$. We will denote by $z_{c_1\ldots c_j}$ such
preimages in $K_i(\Lambda_{(t^\kappa,\Gamma,L^\kappa)})$.

Since $z(0)=0$ we can assume that $(c_1\cdots c_j)_*(z)$ has a
preimage in $K_i(k[{\mathcal N}])/K_i(k)$ for $j\gg0$ with respect
to the right-most arrow in the diagram above.

Fix a grading $k[N_*]=k\oplus R_1\oplus\cdots$ making ellement of
$N_*$ homogenous. It restricts to a grading $k[{\mathcal
N}]=k\oplus S_1\oplus\cdots$. The homomorphism $k[{\mathcal
N}]\xto{\Delta}A\subset M_{2\times2}(k[{\mathcal N}])$ induces a
continuous $\W(k)$-module endomorphism of $K_i(k[{\mathcal
N}])/K_i(k)$ which is a $k$-vector space (see \S\ref{WITT}). On
the  other hand, by Theorem \ref{keat}(a) (and the Berrick-Keating
result on triangular matrices) this endomorphism is just the
multiplication by 2. Therefore, there is a preimage $z^{\mathcal
N}_{c_1\ldots c_j}\in K_i(k[{\mathcal N}])$ of the image of
$z_{c_1\ldots c_j}$ in $K_i(M_{2\times 2}(k[{\mathcal N}]))
=K_i(k[{\mathcal N}])$ under the homomorphism $K_i({\mathcal
N})\xto{\Delta_*}K_i(A)\to K_i(M_{2\times2}(k[{\mathcal N}]))$.

{\it Step 3.} We claim that for $j\gg0$ the preimages
$z_{c_1\ldots c_j}\in K_i(\Lambda_{(t^\kappa,\Gamma,L^\kappa)})$
and $z_{c_1\ldots c_j}^{\mathcal N}\in K_i(k[L])$ can be chosen in
such a way that they map to the same elements in $K_i(A)$.

By Theorem \ref{keat} we have $K_i(A)=K_i(k[{\mathcal N}])\oplus
K_i(k[{\mathcal N}]/\omega^\kappa k[{\mathcal N}])$ and it is
enough to show that the elements $z_{c_1\ldots c_j}$ with $j\gg0$
can be chosen in such a way that their images belong to
$K_i(k[{\mathcal N}])\oplus0$.

First observe, that we can assume $z_{c_1\ldots c_j}(0)=0$. Then,
by Proposition \ref{totaro} and the fact that
$$
(c_1\cdots c_j)_*(z)\in\I{\big(}K_i(k[N_*]_{[c_1\cdots c_j]},
k[N_*]_{[c_1\cdots c_j]}^+)\to K_i(k[N_*],k[N_*]^+){\big)},
$$
we can further choose
$$
z_{c_1\ldots c_j}\in I_{c_1\cdots c_j}(k)
K_i(\Lambda_{(t^\kappa,\Gamma,L^\kappa)},
\Lambda_{(t^\kappa,\Gamma,L^\kappa)}^+)
$$
(w.r.t. to the induced gradings). The image of $z_{c_1\ldots c_j}$
in $K_i(A)$ belongs to the subgroup $I_{c_1\cdots
c_j}(k)K_i(A,A^+)\subset K_i(A,A^+)$. Denote this image by
$(z',z'')$ for the appropriate elements
$$
z'\in I_{c_1\cdots c_j}(k)K_i(k[{\mathcal N}],k[{\mathcal
N}]^+),\quad z''\in I_{c_1\cdots c_j}(k)K_i(k[{\mathcal
N}]/\omega^\kappa k[{\mathcal N}],(k[{\mathcal N}]/\omega^\kappa
k[{\mathcal N}])^+).
$$
Recall, by Theorem \ref{keat} we have the natural identification
$$
K_i(k[{\mathcal N}]/\omega^\kappa k[{\mathcal N}])=
\Ker{\big(}K_i(B)\to K_i(k[{\mathcal N}]/\omega^\kappa k[{\mathcal
N}]){\big)}
$$
where
$$
B=\begin{pmatrix}
k[{\mathcal N}]/\omega^\kappa k[{\mathcal N}]&0\\
\\
k[{\mathcal N}]/\omega^\kappa k[{\mathcal N}]&k[{\mathcal N}]/
\omega^\kappa k[{\mathcal N}]
\end{pmatrix},
$$
and moreover, using this identification, $z''$ is in the image of
the composite map
$$
K_i(\Lambda_{(t^\kappa,\Gamma,L^\kappa)})\to
K_i(\Lambda_{t^\kappa,C})\to K_i(A)\to
K_i(M_{2\times2}(k[{\mathcal N}]))\to K_i(k[{\mathcal
N}]/\omega^\kappa k[{\mathcal N}]).
$$
Now the image of $\Lambda_{(t^\kappa,\Gamma,L^\kappa)}$ in $B$ is
$$
\begin{pmatrix}
\bar\Lambda&0\\
\bar U&\bar\Lambda
\end{pmatrix}
$$
where $\bar\Lambda$ and $\bar U$ denote respectively the images of
$k[L^\kappa(\Gamma)]$ and $k[L^\kappa(\Gamma)]+t^\kappa
k[L^\kappa(\Gamma)]$ in $k[{\mathcal N}]/\omega^\kappa k[{\mathcal
N}]$. Since
$\RR_+\Gamma\setminus\{0\}\subset\int(C_1)\subset\int(C)$ we have
$\dim_k\bar\Lambda<\infty$, the cone $C_1$ being as in Lemma
\ref{bpappro}. (Actually, $\dim_k\bar U<\infty$ as well.) By
Berrick-Keating's result on triangular matrices and Lemma
\ref{findeg} we get the desired vanishing of $z''$ for $j\gg0$.

{\it Step 4.} For any natural number $c$ we have the following
diagram
$$
\xymatrix{
&K_i(\Lambda_{t^\kappa,C})\ar[dd]^{{\tilde c}_*}\ar[rrd]&&&\\
K_i(\Lambda_{(t^\kappa,\Gamma,L^\kappa)})\ar[ru]&&
&K_i(A)\ar@{-->}[r]\ar[dd]^{{\tilde c}_*}&
K_i(M_{2\times2}(k[N_*]))\ar[dd]^{c_*}\\
&K_i(\Lambda_{t^\kappa,C})\ar[rrd]&
K_i(k[{\mathcal N}])\ar[ru]\ar[dd]^{c_*}\ar[rru]&&\\
&&&K_i(A)\ar@{-->}[r]&K_i(M_{2\times2}(k[N_*]))\\
&&K_i(k[{\mathcal N}])\ar[ru]\ar[rru]&&& }
$$
with the obvious homomorphisms. (The endomorphism $\tilde
c:\Lambda_{t^\kappa,C}\to\Lambda_{t^\kappa,C}$ extends naturally
to an endomorphism of $A$, which we denote by the same $\tilde
c$.) All the squares with continuous arrows commute, while there
is no obvious reason why the square with dashed arrows should
commute -- certainly, it does not commute on the level of the
underlying rings. By Step 3 the preimages $z_{c_1\ldots c_j}$ and
$z^{\mathcal N}_{c_1\ldots c_j}$ agree in $K_i(A)$ for $j\gg0$.
Therefore, by Lemma \ref{endo} we achieve the desired pyramidal
descent by fixing $j\gg0$ and running $c$ through $\{c_{j+1}\cdots
c_{j+k}\}_{k\gg0}$. The affine case has been proved.

{\it Step 5.} For the projective varieties we can use the result
in the affine case -- the same arguments as in the proof of
Proposition \ref{global} go through. One only needs to observe
that the standard affine covering of  $\Proj(k[N])$ admits the
following description. It consists of the affine toric varieties
$\Spec(k[N_v])$, defined by the data:
\begin{itemize}
\item
$v$ runs through the vertex set of $\Phi(N)$,
\item
$\lambda$ is a natural number such that $\lambda v\in\gp(N)$,
\item
$C_v\subset\RR\otimes\gp(N)$ is the $(\dim C-1)$-dimensional cone
spanned by $\lambda\Phi(N)$ at $\lambda v$,
\item
$N_v=-\lambda v+(C_v\cap\gp(N))$.
\end{itemize}
The cones $C_v$ are all of type $\BP$ provided $N$ is such,
equivalently -- corner cones of polytopes in $\BP$ are of type
$\BP$. Notice that we also need to take care on monomial
localizations of $k[N]$ -- they also show up in our Mayer-Vietoris
sequences. But such localizations are Laurent polynomial
extensions of monoid rings whose underlying monoids are of type
$\BP$. In particular, we can use the Fundamental Theorem as in the
proof of Proposition \ref{global} (the affine case). \qed

\end{document}